\documentclass[journal]{IEEEtran}
\ifCLASSINFOpdf
\else
\fi
\usepackage[cmex10]{amsmath}
\usepackage[caption=false,font=footnotesize]{subfig}

\usepackage{graphicx}
\usepackage{fixltx2e}
\usepackage{stfloats}
\usepackage{tabularx}

\usepackage{epsfig} 
\usepackage{epstopdf}
\usepackage{amssymb}  

\usepackage{cite}
\usepackage{url}
\usepackage[hidelinks]{hyperref}
\usepackage[hyphenbreaks]{breakurl}
\usepackage{booktabs}
\usepackage[american,fulldiode]{circuitikz} 

\usepackage{algorithm}
\usepackage[noend]{algpseudocode}
\usepackage{float}
\newfloat{algorithm}{t}{lop}

\usepackage{tablefootnote}

\usepackage{tikz}
\usetikzlibrary{calc}

\newcommand{\tikzmark}[1]{\tikz[overlay,remember picture] \node (#1) {};}
\newcommand{\DrawBox}[8][]{%
    \tikz[overlay,remember picture]{%
        \coordinate (TopLeft)     at ($(#2)+(-0.5em,1em)$);
        \coordinate (BottomRight) at ($(#3)+(0.5em,-0.4em)$);
        \path (TopLeft); \pgfgetlastxy{\XCoord}{\IgnoreCoord};
        \path (BottomRight); \pgfgetlastxy{\IgnoreCoord}{\YCoord};
        \draw [#8, #6, rounded corners, draw=#5, fill=#5, fill opacity=#7] (TopLeft) rectangle (BottomRight);
    }
}

\begin{document}
%
\title{Convex Relaxations of Optimal Power Flow Problems: An Illustrative Example}
%
%
%

\author{Daniel~K.~Molzahn,~\IEEEmembership{Member,~IEEE} and 
        Ian~A.~Hiskens,~\IEEEmembership{Fellow,~IEEE}
\thanks{University of Michigan Department of Electrical Engineering and Computer Science: \{molzahn,hiskens\}@umich.edu}
\thanks{Support from the Dow Sustainability Fellowship, ARPA-E grant
  \mbox{DE-AR0000232}, and Los Alamos National Laboratory subcontract
  270958 is gratefully acknowledged.}}

\maketitle

\begin{abstract}
Recently, there has been significant interest in convex relaxations of
the optimal power flow (OPF) problem. A semidefinite programming (SDP) relaxation globally solves many OPF problems. However, there exist practical
problems for which the SDP relaxation fails to yield the
global solution. Conditions for the success or failure of the
SDP relaxation are valuable for determining whether the
relaxation is appropriate for a given OPF problem. To move beyond
existing conditions, which only apply to a limited class of problems,
a typical conjecture is that failure of the SDP relaxation
can be related to physical characteristics of the system. By presenting an example OPF problem with two equivalent formulations, this paper demonstrates that physically based conditions cannot universally explain algorithm behavior. The SDP relaxation fails for one formulation but succeeds in finding the global solution to the other formulation. Since these formulations represent the same system, success (or otherwise) of the SDP relaxation must involve factors beyond just the network physics. The lack of universal physical conditions for success of the SDP relaxation motivates the development of tighter convex relaxations capable of solving a broader class of problems. Tools from polynomial optimization theory provide a means of developing tighter relaxations. We use the example OPF problem to illustrate relaxations from the Lasserre hierarchy for polynomial optimization and a related ``mixed semidefinite/second-order cone programming'' hierarchy.
\end{abstract}


%
\IEEEpeerreviewmaketitle

\section{Introduction}
%
%
%
%
\IEEEPARstart{T}{he} optimal power flow (OPF) problem determines a
minimum cost operating point for an electric power system subject to
both network constraints and engineering limits. Typical objectives
are minimization of losses or generation costs. The OPF problem is
generally non-convex due to the non-linear power flow
equations~\cite{bernie_opfconvexity} and may have local
solutions~\cite{bukhsh_tps}. OPF solution techniques are therefore an
ongoing research topic. Many techniques have been proposed, including
successive quadratic programs, Lagrangian relaxation, and interior
point
methods~\cite{opf_survey,opf_litreview1993IandII,opf_litreview,ferc4,matpower}.

There has been significant interest in convex relaxations of OPF
problems. Convex relaxations lower bound the objective value, can
certify infeasibility, and, in many cases, globally solve OPF
problems. In contrast, traditional OPF solution methods
  may find the global optimum
  \cite{molzahn_lesieutre_demarco-global_optimality_condition} but
  provide no guarantee of doing so, do not provide a measure of
  solution quality and cannot provably identify infeasibility.
Convex relaxations thus have capabilities that supplement traditional
techniques.

For radial systems that satisfy certain non-trivial technical
conditions~\cite{low_tutorial}, a second-order cone programming (SOCP)
relaxation is provably \emph{exact} (i.e., the lower bound is tight
and the solution provides the globally optimal decision
variables). For more general OPF problems, a semidefinite programming
(SDP) based Shor relaxation~\cite{shor-1987b} is often
exact~\cite{lavaei_tps}. Developing tighter and faster relaxations is
an active research area~\cite{sun2015,bienstock2015,coffrin2015}.

Despite success in globally solving many practical OPF
problems~\cite{lavaei_tps,molzahn_holzer_lesieutre_demarco-large_scale_sdp_opf},
there are problems for which the SDP relaxation of~\cite{lavaei_tps}
is not
exact~\cite{allerton2011,bukhsh_tps,hicss2014,molzahn_holzer_lesieutre_demarco-large_scale_sdp_opf}. There
is substantial interest in developing sufficient conditions for
exactness of the SDP relaxation. Existing conditions include
requirements on power injection, voltage magnitude and line-flow
limits, and either radial networks (typical of distribution systems),
appropriate placement of controllable phase-shifting transformers, or
a limited subset of mesh network
topologies~\cite{low_tutorial,madani2014}.

The SDP relaxation globally solves many OPF problems which do not
satisfy any known sufficient
conditions~\cite{low_tutorial,madani2014}. In other words, the set of
problems guaranteed to be exact by known sufficient conditions is much
smaller than the actual set of problems for which the relaxation is
exact. This suggests the potential for developing less
  stringent conditions. It is natural to speculate that some physical
characteristics of an OPF problem may govern such conditions. With
solutions that are close to voltage collapse, this speculation is
supported by several problems for which the SDP relaxation fails to be
exact~\cite{hicss2014}.

This paper presents an example that dampens enthusiasm for this avenue of research. We consider a small problem, first presented in~~\cite{iscas2015}, with two equivalent formulations. The SDP relaxation globally solves one formulation but fails to solve the other. Since both formulations represent the same system, strictly physically based conditions for the success of the relaxation cannot differentiate between these formulations.\footnote{See also~\cite{madani2014} for an example where different line-flow limit formulations determine success or failure of the SDP relaxation.} The feasible spaces of both problems illustrate why the SDP relaxation succeeds for one formulation but fails for the other.

The small example considered in this paper is relatively simple. In fact, this example ``OPF'' problem reduces to finding the minimum loss solution to power flow constraint equations for a specified set of power injections and voltage magnitudes. Thus, this example further demonstrates that the SDP relaxation may fail even for simple OPF problems.

The lack of universal, physically based conditions for determining success or failure of the SDP relaxation of~\cite{lavaei_tps} motivates the development of tighter convex relaxations. Recent research~\cite{pscc2014,cedric,ibm_paper,molzahn_hiskens-sparse_moment_opf,cdc2015} exploits the fact that the OPF problem is a polynomial optimization problem in terms of the complex voltage phasors. Separating the complex voltages into real and imaginary parts yields a polynomial optimization problem in real variables. This facilitates the application of the Lasserre hierarchy of ``moment'' relaxations for real polynomial optimization problems, which take the from of SDPs~\cite{lasserre_book}. The first-order moment relaxation is equivalent to the SDP relaxation of~\cite{lavaei_tps}. Higher-order moment relaxations thus generalize the SDP relaxation of~\cite{lavaei_tps}.

Moment relaxations globally solve many problems for which the SDP relaxation of~\cite{lavaei_tps} is not exact~\cite{pscc2014,cedric,ibm_paper,molzahn_hiskens-sparse_moment_opf,cdc2015}. We use the small example system to illustrate the moment relaxations, including exploration of the feasible space.

The ability of the moment relaxations to solve a broader class of OPF
problems comes at a computational cost: the matrices grow rapidly with
both relaxation order and system size. Ameliorating the former
challenge, low relaxation orders suffice for global solution of many
problems. Several recent developments address the latter
challenge. First, by exploiting sparsity and selectively applying the
higher-order constraints to specific buses, loss-minimization problems
with thousands of buses are computationally
tractable~\cite{molzahn_hiskens-sparse_moment_opf,cdc2015}. Second,
rather than separating complex voltages into their real and imaginary
parts, a hierarchy built directly from the complex formulation is
computationally advantageous~\cite{complex_hierarchy}. Third, a
``mixed SDP/SOCP'' hierarchy implements the first-order constraints
with a SDP formulation, but the higher-order constraints are relaxed
to a SOCP formulation. The less computationally intensive SOCP
constraints often reduce solution times while still yielding global
optima. The small example system is used to illustrate the mixed
SDP/SOCP hierarchy.

This paper is organized as follows. Section~\ref{l:opf_formulation} introduces the OPF problem. Section~\ref{l:sdp_relaxation} describes the SDP relaxation of~\cite{lavaei_tps}. Section~\ref{l:example} presents the example OPF problem which demonstrates that factors beyond the problem physics determine success or failure of the SDP relaxation. Sections~\ref{l:moment} and~\ref{l:mixed_sdpsocp} provide the moment relaxations using SDP constraints and the mixed SDP/SOCP hierarchy, respectively, with the problem in Section~\ref{l:example} providing an illustrative example. Section~\ref{l:conclusion} concludes the paper.

\section{Optimal Power Flow Problem}
\label{l:opf_formulation}

We first present an OPF formulation in terms of rectangular voltage
coordinates, active and reactive power injections, and apparent power
line flow limits. Consider an $n$-bus system with $n_l$ lines, where
$\mathcal{N} = \left\lbrace 1, \ldots, n \right\rbrace$ is the set of
buses, $\mathcal{G}$ is the set of generator buses, and $\mathcal{L}$
is the set of lines. The network admittance matrix is $\mathbf{Y} =
\mathbf{G} + \mathbf{j} \mathbf{B}$, where $\mathbf{j}$ denotes the
imaginary unit. Let $P_{Dk} + \mathbf{j} Q_{Dk}$ represent the active
and reactive load demand and $V_k = V_{dk} + \mathbf{j} V_{qk}$ the
voltage phasors at each bus~$k \in \mathcal{N}$. Superscripts ``max''
and ``min'' denote specified upper and lower limits. Buses without
generators have maximum and minimum generation set to zero. 


Define a function for squared voltage magnitude:
\begin{equation} \label{opf_Vsq}
f_{Vk}\left(V_d, V_q\right) := V_{dk}^2 + V_{qk}^2.
\end{equation}
The power flow equations describe the network physics:
\begin{subequations}
\label{opf_balance}
\begin{align}\nonumber
f_{Pk}\left(V_d,V_q\right) := & V_{dk} \sum_{i=1}^n \left( \mathbf{G}_{ki} V_{di} - \mathbf{B}_{ki} V_{qi} \right) &  &  \\[-5pt] 
\label{opf_Pbalance}  & + V_{qk} \sum_{i=1}^n \left( \mathbf{B}_{ki}V_{di} + \mathbf{G}_{ki}V_{qi} \right) + P_{Dk}, \\ \nonumber 
f_{Qk}\left(V_d,V_q\right) := & V_{dk} \sum_{i=1}^n \left( -\mathbf{B}_{ki}V_{di} - \mathbf{G}_{ki} V_{qi}\right) \\[-5pt]
\label{opf_Qbalance} & + V_{qk} \sum_{i=1}^n \left( \mathbf{G}_{ki} V_{di} - \mathbf{B}_{ki} V_{qi}\right) + Q_{Dk}.
\end{align}
\end{subequations}
Define a convex quadratic cost of active power generation:
\begin{equation}\label{objfunction}
f_{Ck}\left(V_d,V_q\right) := c_{k2} \left(f_{Pk}\left(V_d,V_q\right)\right)^2 + c_{k1} f_{Pk}\left(V_d,V_q\right) + c_{k0}.
\end{equation}

\begin{figure}[t]
\centering
\begin{circuitikz}[scale=0.85]\draw
  (0,2) to[short,o-] (1,2)
  		to[cute inductor] (1,0)
  		to[short,-o](0,0)
  (0,2) to[open,v_=$V_l$] (0,0)
  (3,0) to[short] (2,0)
  		to[cute inductor] (2,2)
  		to[short] (3,2)
  (1.5,0) node[anchor=north]{$\tau_{lm} e^{\mathbf{j}\theta_{lm}}:1$}
  (2,0) node[anchor=east] {}
        to[short] (7,0)
  (2,2) node[anchor=east] {}
  		to[short] (3,2)
  		to[R=$R_{lm}$] (5,2)
  		to[cute inductor,l=$\mathbf{j}X_{lm}$] (7,2)
  		to[C,l_=$\mathbf{j}\frac{b_{sh,lm}}{2}$] (7,0)
  (3,0) to[C,l_=$\mathbf{j}\frac{b_{sh,lm}}{2}$] (3,2)
  (7,2) to[short,-o] (8,2)
  (7,0) to[short,-o] (8,0)
  (8,2) to[open,v^=$V_m$] (8,0)
  (0,2.3) node[anchor=south]{$P_{lm} + \mathbf{j} Q_{lm}$}
  (8,2.3) node[anchor=south]{$P_{ml} + \mathbf{j} Q_{ml}$};
  \draw[thick,->] (8.1,2.3) -- (7,2.3);
  \draw[thick,->] (-0.25,2.3) -- (0.75,2.3);
\end{circuitikz}
\vspace{-10pt}
\caption{Line Model} \label{f:linemodel}
\end{figure}
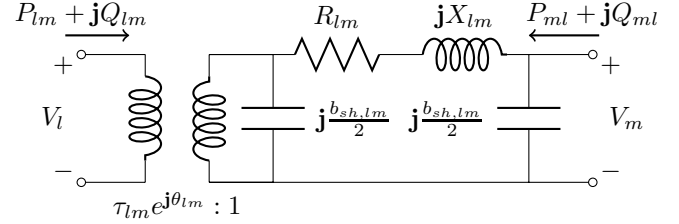

We use a line model with an ideal transformer that has a specified turns ratio $\tau_{lm} e^{\mathbf{j}\theta_{lm}} \colon 1$ in series with a $\Pi$ circuit with series impedance $R_{lm} + \mathbf{j} X_{lm}$ (equivalent to an admittance of $g_{lm} + \mathbf{j} b_{lm} := \frac{1}{R_{lm} + \mathbf{j} X_{lm}}$) and total shunt susceptance $\mathbf{j} b_{sh,lm}$. (See Fig.~\ref{f:linemodel}.) The line flow equations are:
\begin{subequations}
\begin{align}
\nonumber & f_{Plm}\left(V_d,V_q\right) := \left( V_{dl}^2 + V_{ql}^2\right) g_{lm}/\tau_{lm}^2 \\ \nonumber
& \; +
\left(V_{dl}V_{dm} + V_{ql}V_{qm}\right) \left(b_{lm}\sin\left(\theta_{lm} \right) - g_{lm}\cos\left(\theta_{lm} \right) \right) / \tau_{lm} \\
\label{Plm}& \; +
\left(V_{dl}V_{qm} - V_{ql}V_{dm}\right) \left(g_{lm}\sin\left(\theta_{lm}\right) + b_{lm}\cos\left(\theta_{lm}\right)\right) / \tau_{lm},  \\
\nonumber & f_{Qlm}\left(V_d,V_q\right) := -\left( V_{dl}^2 + V_{ql}^2\right) \left(b_{lm} + \frac{b_{sh,lm}}{2}\right) / \tau_{lm}^2 \\ \nonumber 
& \; +
\left(V_{dl}V_{dm} + V_{ql}V_{qm}\right) \left(b_{lm}\cos\left(\theta_{lm} \right) + g_{lm}\sin\left(\theta_{lm} \right) \right) / \tau_{lm} \\
\label{Qlm} & \; +
\left(V_{dl}V_{qm} - V_{ql}V_{dm}\right) \left(g_{lm}\cos\left(\theta_{lm}\right) - b_{lm}\sin\left(\theta_{lm}\right)\right) / \tau_{lm}, \\
\label{Slm} & f_{Slm}\left(V_d, V_q \right) := \left(f_{Plm}\left(V_d, V_q \right)\right)^2 + \left(f_{Qlm}\left(V_d, V_q \right)\right)^2,\\
\nonumber & f_{Pml}\left(V_d,V_q\right) := \left(V_{dm}^2 + V_{qm}^2 \right)g_{lm} \\ \nonumber & \; -
\left(V_{dl}V_{dm} + V_{ql}V_{qm} \right)\left(g_{lm}\cos\left(\theta_{lm}\right) + b_{lm}\sin\left(\theta_{lm}\right) \right) / \tau_{lm} \\ \label{Pml}
& \; + 
\left(V_{dl}V_{qm} - V_{ql}V_{dm} \right) \left(g_{lm}\sin\left(\theta_{lm}\right) - b_{lm}\cos\left(\theta_{lm}\right) \right) / \tau_{lm},
\\ 
\nonumber & f_{Qml}\left(V_d,V_q\right) := -\left( V_{dm}^2 + V_{qm}^2\right) \left(b_{lm} + \frac{b_{sh,lm}}{2}\right) \\ \nonumber & \; +
\left(V_{dl}V_{dm} + V_{ql}V_{qm}\right) \left(b_{lm}\cos\left(\theta_{lm} \right) - g_{lm}\sin\left(\theta_{lm} \right) \right)  / \tau_{lm} \\ \label{Qml} & \; + 
\left(-V_{dl}V_{qm} + V_{ql}V_{dm}\right) \left(g_{lm}\cos\left(\theta_{lm}\right) + b_{lm}\sin\left(\theta_{lm}\right)\right) / \tau_{lm}, \\
\label{Sml} & f_{Sml}\left(V_d, V_q \right) := \left(f_{Pml}\left(V_d, V_q \right)\right)^2 + \left(f_{Qml}\left(V_d, V_q \right)\right)^2.
\end{align}
\end{subequations}

The OPF problem is:
\begin{subequations}
\label{opf}
\begin{align}
\label{opf_obj} & \min_{V_d,V_q}\quad \sum_{k \in \mathcal{G}} f_{Ck}\left(V_d, V_q \right) \qquad \mathrm{subject\; to} \hspace{-160pt} & \\
\label{opf_P} &  \quad P_{Gk}^{\mathrm{min}} \leqslant f_{Pk}\left(V_d,V_q\right) \leqslant P_{Gk}^{\mathrm{max}} & \forall k \in \mathcal{N} \\
\label{opf_Q} &  \quad Q_{Gk}^{\mathrm{min}} \leqslant f_{Qk}\left(V_d,V_q\right) \leqslant Q_{Gk}^{\mathrm{max}} &  \forall k \in \mathcal{N} \\
\label{opf_V} &  \quad (V_{k}^{\mathrm{min}})^2 \leqslant f_{Vk}\left(V_d,V_q\right) \leqslant (V_{k}^{\mathrm{max}})^2 &  \forall k \in \mathcal{N}  \\
\label{opf_Slm} & \quad f_{Slm}\left(V_d,V_q\right) \leqslant \left(S_{lm}^{\mathrm{max}}\right)^2 &  \forall \left(l,m\right) \in \mathcal{L} \\ 
\label{opf_Sml} & \quad f_{Sml}\left(V_d,V_q\right) \leqslant \left(S_{lm}^{\mathrm{max}}\right)^2 &  \forall \left(l,m\right) \in \mathcal{L} \\ 
\label{opf_Vref} & \quad V_{q1} = 0.
\end{align}
\end{subequations}
Constraint~\eqref{opf_Vref} sets the reference bus angle to zero.

\section{Semidefinite Relaxation of the OPF Problem}
\label{l:sdp_relaxation}

This section describes a SDP relaxation of the OPF problem adopted
from~\cite{lavaei_tps,molzahn_holzer_lesieutre_demarco-large_scale_sdp_opf,andersen2014}. We
use notation from~\cite{pscc2014,molzahn_hiskens-sparse_moment_opf}
corresponding to the moment relaxations that will be introduced in the
following sections. We begin with several definitions. Define the
vector of real decision variables $\hat{x} \in \mathbb{R}^{2n}$ as
\begin{equation}\label{xhat}
\hat{x} := \begin{bmatrix} V_{d1} & V_{d2} & \ldots V_{qn} \end{bmatrix}^\intercal
\end{equation}
\noindent where $\left(\cdot\right)^\intercal$ denotes the transpose.\footnote{The ability to arbitrarily set an angle reference in the OPF problem enables the choice of one arbitrarily selected variable. We choose $V_{q1} = 0$ as in~\eqref{opf_Vref}.} A monomial is defined using a vector $\alpha \in \mathbb{N}^{2n}$ of exponents: $\hat{x}^\alpha := V_{d1}^{\alpha_1} V_{d2}^{\alpha_2}\cdots V_{qn}^{\alpha_{2n}}$. A polynomial is $g\left(\hat{x}\right) := \sum_{\alpha \in \mathbb{N}^{2n}} g_{\alpha} \hat{x}^{\alpha}$, where $g_{\alpha}$ is the real scalar coefficient corresponding to the monomial $\hat{x}^{\alpha}$.

Define a linear functional $L_y\left\lbrace g\right\rbrace$ which replaces the monomials $\hat{x}^{\alpha}$ in a polynomial $g\left(\hat{x}\right)$ with real scalar variables $y$:
\begin{equation}
\label{eq:Lreal}
L_y\left\lbrace g \right\rbrace := \sum_{\alpha \in \mathbb{N}^{2n}} g_{\alpha} y_{\alpha}.
\end{equation}
\noindent For a matrix $g\left(\hat{x}\right)$, $L_y\left\lbrace g\right\rbrace$ is applied componentwise to each element of $g\left(\hat{x}\right)$. 

Consider, for example, the vector $\hat{x} = \begin{bmatrix}V_{d1} & V_{d2} & V_{q2} \end{bmatrix}^\intercal$ corresponding to the voltage components of a two-bus system, where the angle reference~\eqref{opf_Vref} is used to eliminate $V_{q1}$. Consider also the polynomial $g\left(\hat{x}\right) = \left(V_2^{\max}\right)^2 - V_{d2}^2 - V_{q2}^2$. (The constraint $g\left(\hat{x}\right) \geqslant 0$ forces the voltage magnitude at bus~2 to be less than or equal to $V_2^{\max}$.) Then $L_y\left\lbrace g\right\rbrace = \left(V_2^{\max}\right)^2y_{000} - y_{020} - y_{002}$. Thus, $L_y\left\lbrace g \right\rbrace$ converts a polynomial $g\left(\hat{x}\right)$ to a linear function of $y$.



The SDP relaxation of~\eqref{opf} is:
\begin{subequations}
\label{sdpprimal}
\begin{align}
\label{sdp_obj} & \min_{y,\omega} \sum_{k \in \mathcal{G}} \omega_k \qquad \mathrm{subject\; to} \hspace{-35pt} & \\
\label{sdp_P} &  \quad P_{Gk}^{\mathrm{min}} \leqslant L_y\left\lbrace f_{Pk} \right\rbrace \leqslant P_{Gk}^{\mathrm{max}} & \forall k \in \mathcal{N} \\
\label{sdp_Q} &  \quad Q_{Gk}^{\mathrm{min}} \leqslant L_y\left\lbrace f_{Qk} \right\rbrace \leqslant Q_{Gk}^{\mathrm{max}} &  \forall k \in \mathcal{N} \\
\label{sdp_V} &  \quad \left(V_{k}^{\mathrm{min}}\right)^2 \leqslant L_y\left\lbrace f_{Vk} \right\rbrace \leqslant \left(\vphantom{V_{k}^{\mathrm{min}}} V_{k}^{\mathrm{max}}\right)^2 &  \forall k \in \mathcal{N}  \\
\nonumber & \quad \left(1-c_{k1}L_y\left\lbrace f_{Pk} \right\rbrace-c_{k0} + \omega_k \right) \\ \label{sdp_quadcost} & \qquad \geqslant \left|\left| \begin{bmatrix} \left(1+c_{k1}L_y\left\lbrace f_{Pk} \right\rbrace+c_{k0}-\omega_k \right) \\ 2\sqrt{c_{k2}}\, L_y\left\lbrace f_{Pk} \right\rbrace \end{bmatrix} \right|\right|_2 \hspace{-8pt} & \forall k \in \mathcal{G} \end{align}\begin{align}
\label{sdp_Slm} & \quad S_{lm}^\mathrm{max} \geqslant \left|\left| \begin{bmatrix} L_y\left\lbrace f_{Plm} \right\rbrace \\ L_y\left\lbrace f_{Qlm} \right\rbrace  \end{bmatrix} \right|\right|_2 & \forall \left(l,m \right) \in \mathcal{L} \\
\label{sdp_Sml} & \quad S_{lm}^\mathrm{max} \geqslant \left|\left| \begin{bmatrix} L_y\left\lbrace f_{Pml} \right\rbrace \\ L_y\left\lbrace f_{Qml} \right\rbrace \end{bmatrix} \right|\right|_2 & \forall \left(l,m \right) \in \mathcal{L} \\
\label{sdp_W} & \quad L_y\left\lbrace \hat{x} \hat{x}^\intercal \right\rbrace \succcurlyeq 0 \\
\label{sdp_y0} & \quad y_{00\ldots 0} = 1 & \\
\label{sdp_Vref} & \quad y_{\star\star\ldots\star\rho\star\ldots\star} = 0 & \rho = 1,2,
\end{align}
\end{subequations}
\noindent where $\succcurlyeq 0$ indicates positive semidefiniteness of the corresponding matrix and $\left|\left|\,\cdot\,\right|\right|_2$ denotes the two-norm. The generation cost function \eqref{opf_obj} for the generator at bus~$k$ is implemented using the variable $\omega_k$ and the SOCP formulation in~\eqref{sdp_quadcost}~\cite{andersen2014}. The apparent power line flow constraints~\eqref{opf_Slm} and~\eqref{opf_Sml} are implemented with the SOCP formulations in~\eqref{sdp_Slm} and \eqref{sdp_Sml}. See~\cite{molzahn_holzer_lesieutre_demarco-large_scale_sdp_opf} for a more general formulation of the SDP relaxation that considers the possibilities of multiple generators per bus and convex piecewise-linear generation costs. The constraint~\eqref{sdp_y0} enforces the fact that $\hat{x}^{0} = 1$. The constraint~\eqref{sdp_Vref} corresponds to the angle reference $V_{q1} = 0$; the $\rho$ in \eqref{sdp_Vref} is in the index $n+1$, which corresponds to the variable $V_{q1}$. Note that the angle reference can alternatively be used to eliminate all terms corresponding to $V_{q1}$ to reduce the size of the semidefinite program.

If the condition $\mathrm{rank}\left(L_y\left\lbrace \hat{x} \hat{x}^\intercal \right\rbrace\right) = 1$ is satisfied, the relaxation is ``exact'' and the global solution to~\eqref{opf} is recovered using an eigen-decomposition. Consider a solution to~\eqref{sdpprimal} where the rank of the matrix $L_y\left\lbrace \hat{x} \hat{x}^\intercal \right\rbrace$ is equal to one with non-zero eigenvalue $\lambda$ and associated unit-length eigenvector $\eta$. The globally optimal voltage phasor solution to~\eqref{opf} is
\begin{equation}\label{Vstar}
V^\ast = \sqrt{\lambda} \left(\eta_{1:n} + \mathbf{j} \eta_{\left(n+1\right):2n} \right)
\end{equation}
\noindent where subscripts denote vector entries in MATLAB notation. 

The computational bottleneck of the SDP relaxation is the constraint~\eqref{sdp_W}, which enforces positive semidefiniteness of a $2n \times 2n$ matrix. Solving the SDP relaxation of large OPF problems requires exploiting network sparsity. A matrix completion decomposition exploits sparsity by converting the positive semidefinite constraint on the large matrix in~\eqref{sdp_W} to positive semidefinite constraints on many smaller submatrices. These submatrices are defined using the cliques (i.e., completely connected subgraphs) of a chordal extension of the power system network graph. See~\cite{jabr2011,molzahn_holzer_lesieutre_demarco-large_scale_sdp_opf,andersen2014} for a full description of a formulation that enables solution of~\eqref{sdpprimal} for systems with thousands of buses.

\section{Equivalent Formulations of a Small Example Problem}
\label{l:example}

The SDP relaxation~\eqref{sdpprimal} globally solves many OPF problems which do not satisfy any known sufficient conditions guaranteeing exactness~\cite{low_tutorial,madani2014}, indicating the potential for development of broader sufficient conditions. One speculation is that some physical characteristic of the OPF problem predicts the relaxation's success or failure.

The following example shows that strictly physically based sufficient
conditions are unable to definitively predict success or failure of
the SDP relaxation for all OPF problems. The example problem has
equivalent two- and three-bus formulations. The relaxation globally
solves the two-bus system. For the three-bus system, however, the relaxation only
gives a strict lower bound on the objective value rather than the
solution.

\subsection{Example Problem}

\begin{figure}[t]
\centering
\begin{tikzpicture}[scale=0.9, transform shape]

\path[draw,line width=4pt] (0,0) -- (0,2);
\draw (0,0.4) node[below right] {$1$};
\draw (-1.6,2.3) node[right] {$V_1 = 1$};
\draw (-1.6,1.8) node[right] {$\theta_1 = 0^\circ$};
\path[draw,line width=2pt] (-0.4,1) -- (0,1);
\draw[line width=1] (-0.8,1) circle (0.4);

\path[draw,line width=2pt] (0,1.3) -- (5,1.3);
\draw (1.6,2.4) node[below] {$R_{12} + \mathbf{j} X_{12}$};
\draw (2.5,1.8) node[below] {$= 0.15 + \mathbf{j} 0.1$};

\path[draw,line width=4pt] (5,0) -- (5,2);
\draw (5,0.4) node[below right] {$2$};
\path[draw,line width=2pt] (5,1) -- (5.4,1);
\draw[line width=1] (5.8,1) circle (0.4);
\draw (5.2,2.3) node[right] {$V_2 = 1.3$};
\draw (5.2,1.8) node[right] {$P_2 = 0$};

\path[draw,line width=2pt] (0,0.65) -- (0.5,0.65) -- (1.75,-1.5) -- (1.75,-2);
\draw (-0.2,-0.4) node[below] {$R_{13} + \mathbf{j} X_{13}$};
\draw (0.3,-0.8) node[below] {$= 0.1 + \mathbf{j} 0.05$};

\path[draw,line width=4pt] (1.5,-2) -- (3.5,-2);
\draw (1.5,-2) node[left] {$3$};
\draw[->,line width=2pt] (2.5,-2) -- (2.5,-3);
\draw[line width=1] (5.8,1) circle (0.4);
\draw (2.7,-2.3) node[right] {$P_3 = 0$};
\draw (2.7,-2.8) node[right] {$Q_3 = 0$};

\path[draw,line width=2pt] (5,0.65) -- (4.5,0.65) -- (3.25,-1.5) -- (3.25,-2);
\draw (5,-0.4) node[below] {$R_{23} + \mathbf{j} X_{23}$};
\draw (5.3,-0.8) node[below] {$= 0.001 + \mathbf{j} 0.05$};

\end{tikzpicture}
\caption{Three-Bus System}
\label{f:threebussystem}

\end{figure}
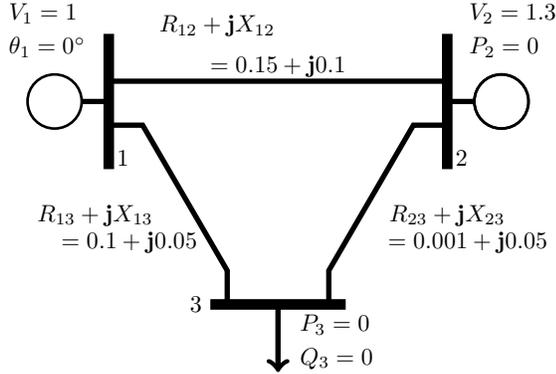

Consider the two- and three-bus systems in Figs.~\ref{f:twobussystem}
and~\ref{f:threebussystem}. For both systems, the voltage magnitudes
at buses~1 and~2 are fixed to $1.0$ and $1.3$~per unit, respectively,
and the active power injection at bus~2 is fixed to
zero.\footnote{Equality constraints are achieved by setting
  the upper and lower limits equal (e.g., $V_1^{\max} = V_1^{\min} =
  1$~per unit).} There are no limits on the reactive power injections
at buses~1 and~2. For bus~3 in the three-bus system, the active and
reactive power injections are constrained to zero and there is no
voltage magnitude constraint. With the active power injections at the
other buses fixed to zero, the objective function minimizes active
power injection at bus~1.

The resistance-to-reactance ratios for lines in both the two- and three-bus systems are somewhat atypical for transmission systems, but are not particularly unusual for more lossy networks like subtransmission and distribution systems~\cite{willis2004}. Similar characteristics to these systems may also occur when using ``equivalencing'' techniques to reduce larger systems to a smaller representative network~\cite{deckmann1980studies,deckmann1980numerical}.

\begin{figure}[t]
\centering
\begin{tikzpicture}[scale=0.9, transform shape]

\path[draw,line width=4pt] (0,0) -- (0,2);
\draw (0,0.4) node[below right] {$1$};
\draw (-1.6,2.3) node[right] {$V_1 = 1$};
\draw (-1.6,1.8) node[right] {$\theta_1 = 0^\circ$};
\path[draw,line width=2pt] (-0.4,1) -- (0,1);
\draw[line width=1] (-0.8,1) circle (0.4);

\path[draw,line width=2pt] (0,1) -- (5,1);
\draw (1.3,2.4) node[below] {$R_{12}^\prime + \mathbf{j} X_{12}^\prime$};
\draw (2.5,1.8) node[below] {$= 0.06129 + \mathbf{j} 0.05117$};

\path[draw,line width=4pt] (5,0) -- (5,2);
\draw (5,0.4) node[below right] {$2$};
\path[draw,line width=2pt] (5,1) -- (5.4,1);
\draw[line width=1] (5.8,1) circle (0.4);
\draw (5.2,2.3) node[right] {$V_2 = 1.3$};
\draw (5.2,1.8) node[right] {$P_2 = 0$};

\end{tikzpicture}
\caption{Two-Bus System}
\label{f:twobussystem}

\end{figure}
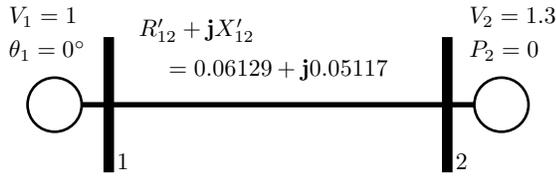

%
%
%
%
%
%

With two quantities specified at each bus~$k$ along with two degrees
of freedom ($V_{dk}$ and $V_{qk}$), the feasible space for the OPF problem~\eqref{opf} for this example consists of a set of isolated points that are the solutions of the power flow equations. The OPF finds the solution point that has the lowest active power losses. Here, this solution corresponds to the ``high-voltage/small angle-difference'' power flow solution, which is commonly calculated using a Newton-Raphson iteration initialized from a flat start (i.e., voltages of $1\angle 0^\circ$).\footnote{For both two- and three-bus systems, \eqref{opf} has one other local minimum: there exists one ``low-voltage/large angle-difference'' power flow solution with larger losses.} In this paper, however, we use this problem to explore the properties of the convex relaxations.

\begin{figure*}[t]
\centering
\subfloat[Projection of the Two-Bus System's Feasible Space]{\includegraphics[totalheight=0.23\textheight]{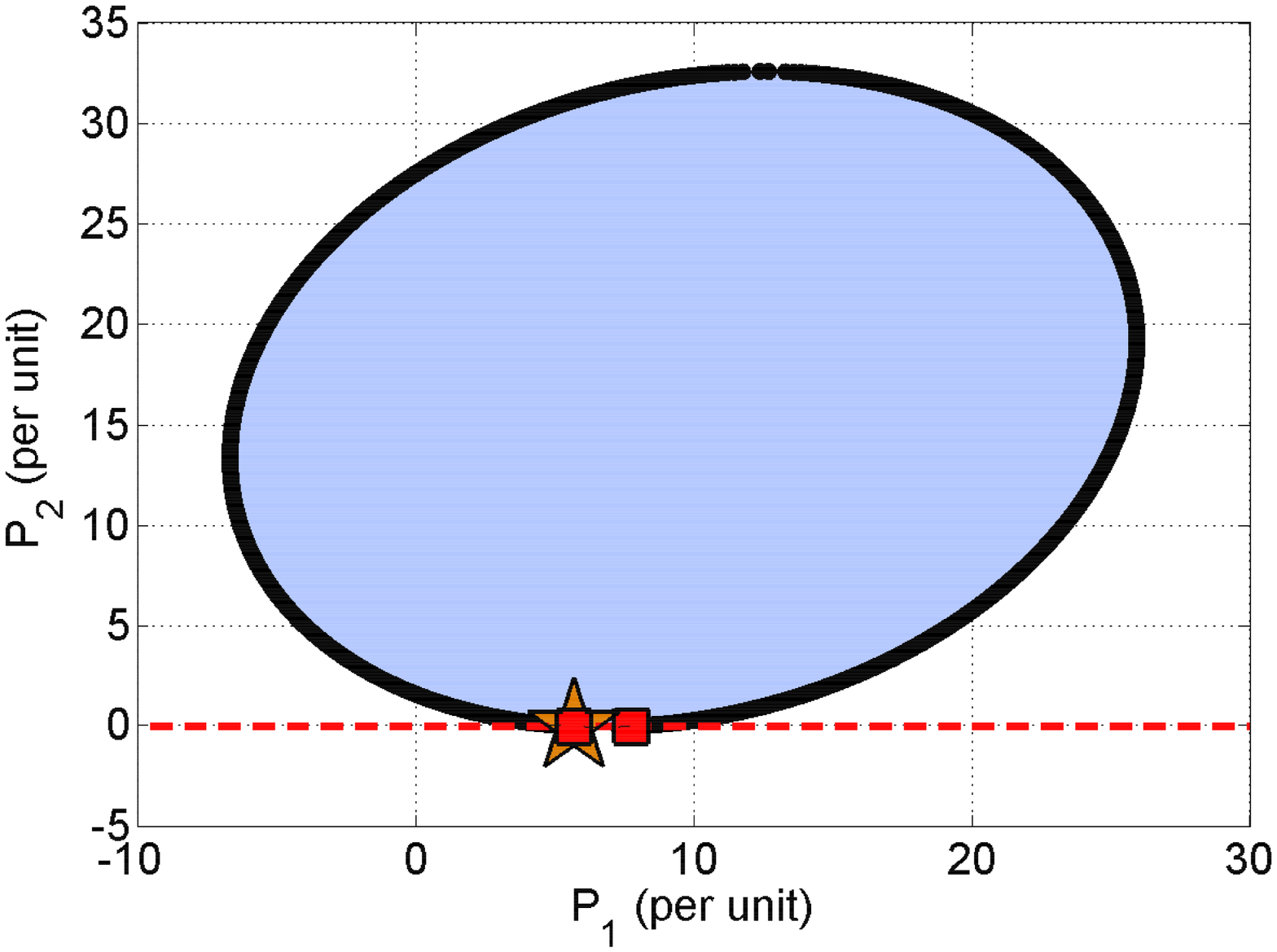}\label{f:twobus_fullspace}}
{\captionsetup{margin={2cm,0cm}}\hspace{-5pt}\subfloat[Zoomed View of Fig.~\ref{f:twobus_fullspace}]{\includegraphics[trim=-5cm 0cm 0cm 0cm, clip=true,totalheight=0.24\textheight]{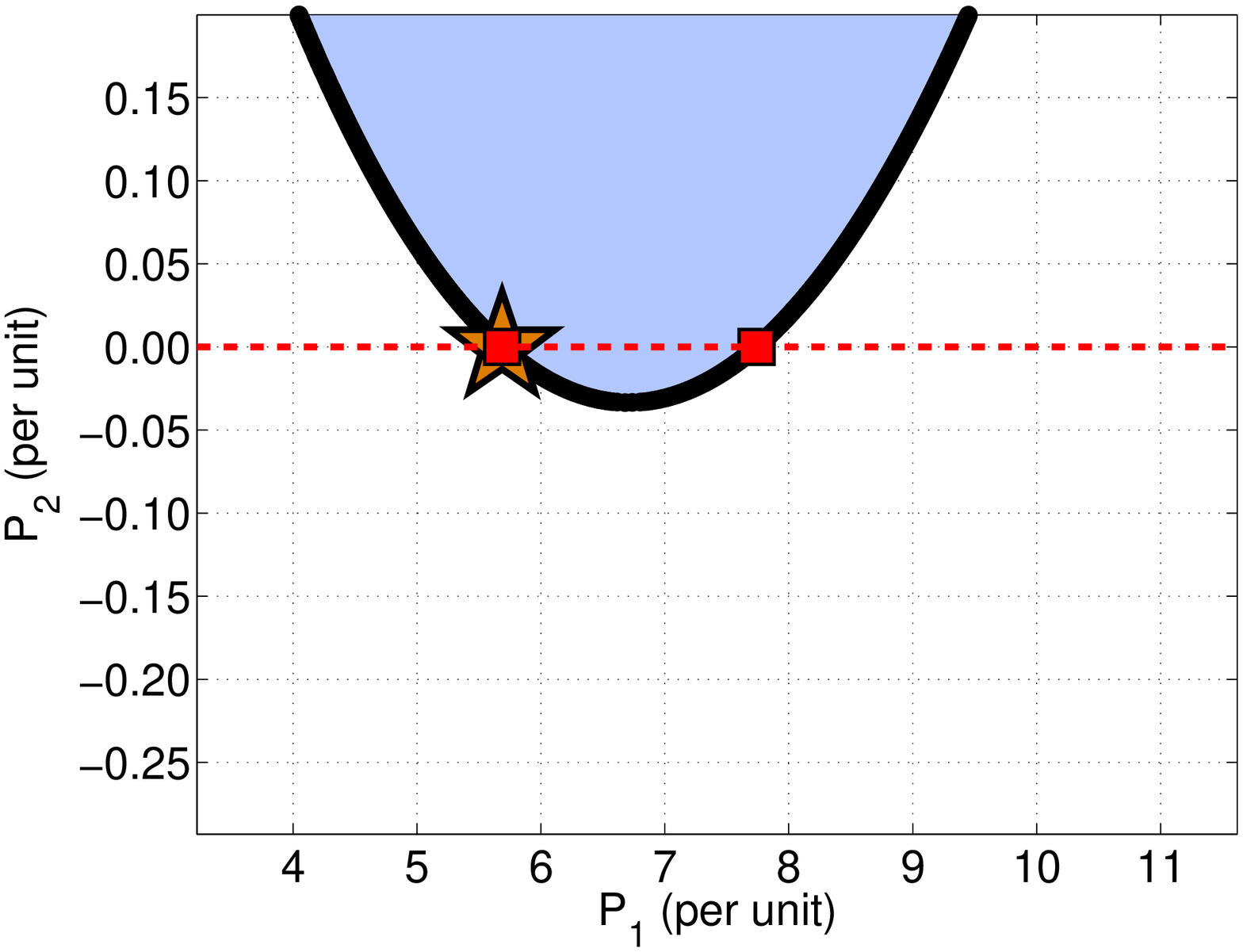}\label{f:twobus_zoomspace}}} \\
\caption{Projection of the Two-Bus System's Feasible Space. The red
  squares at the intersection of the black oval and red dashed line
  are the feasible space for the OPF problem~\eqref{opf}. The blue
  region, including the black oval boundary, is the
  feasible space for the SDP relaxation~\eqref{sdpprimal}. The orange
  star is the solution to the SDP relaxation, which is the global
  optimum for the two-bus system.}
\label{f:twobus_space}
\end{figure*}

\begin{figure*}[t]
\centering 
\subfloat[Projection of the Three-Bus System's Feasible Space]{\includegraphics[trim=0.cm 0cm 0.2cm 3.5cm, clip=true,totalheight=0.20\textheight]{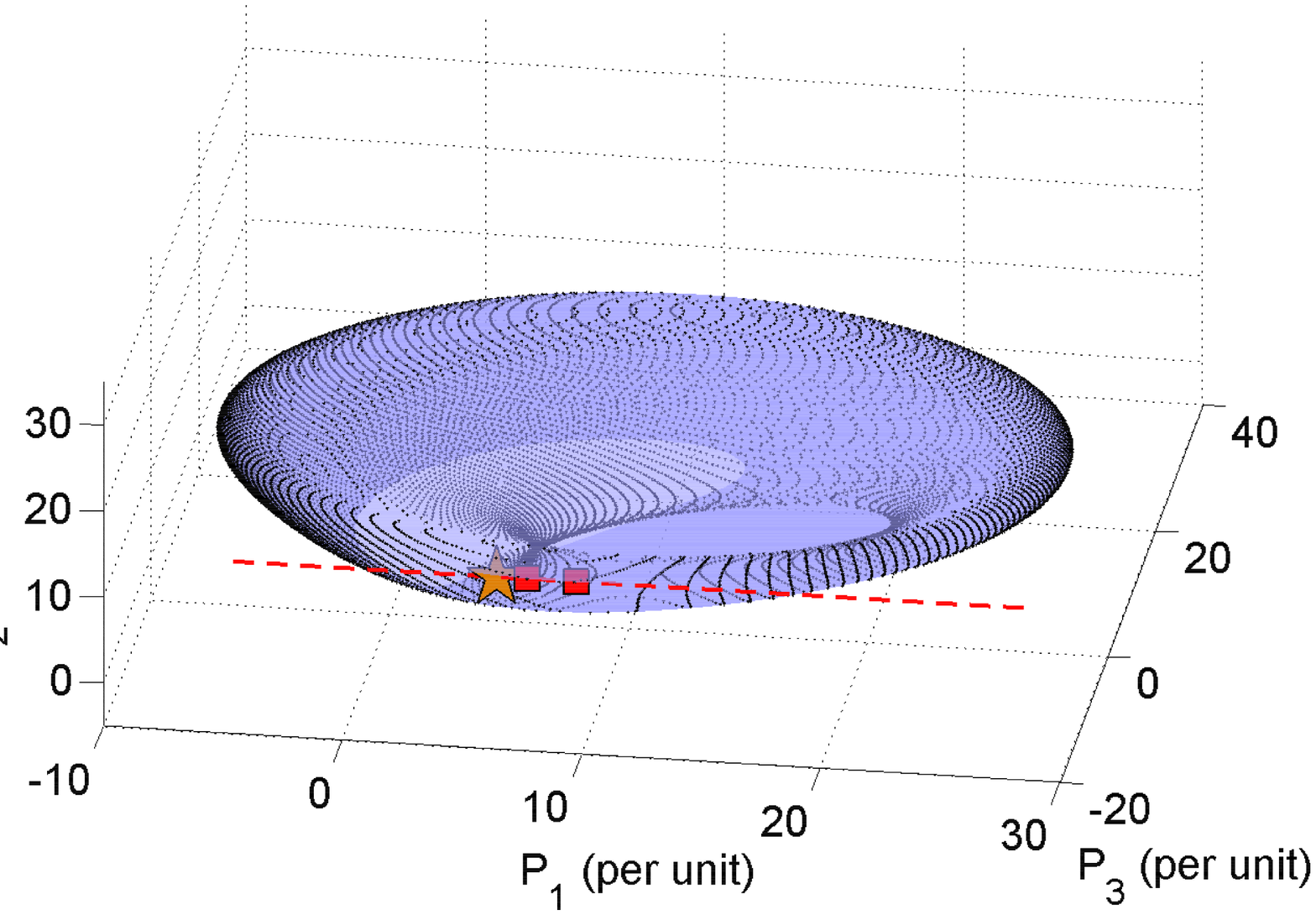}\label{f:threebus_fullspace}} 
\hspace{0pt}
\subfloat[Projection of the Three-Bus System's Feasible Space with $P_{3} = 0$]{\includegraphics[totalheight=0.23\textheight]{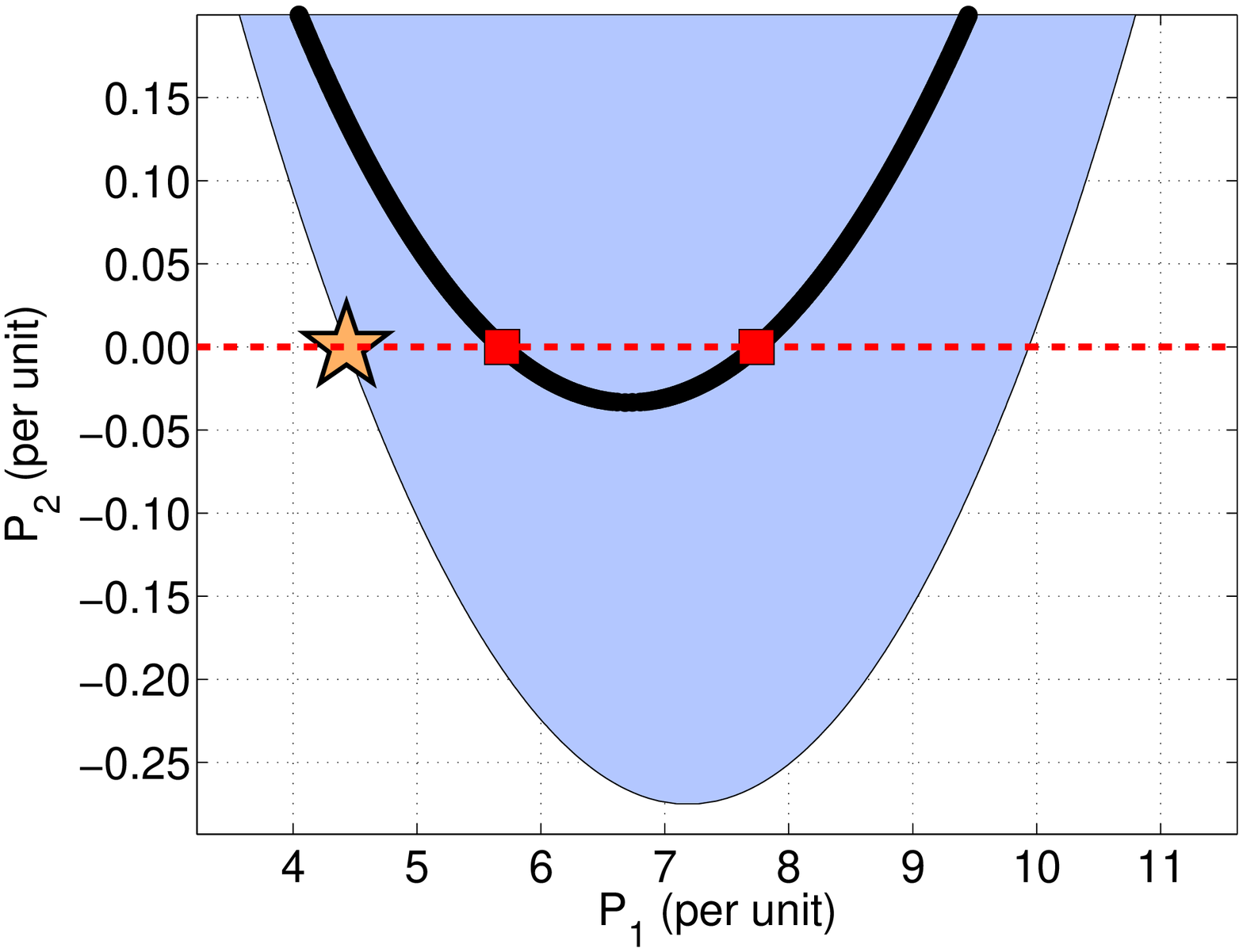}\label{f:threebus_zoomspace}} \\
\caption{Projection of the Three-Bus System's Feasible Space. The
  feasible space for the OPF problem~\eqref{opf} is denoted by the red
  squares at the intersection of the red dashed line and the region
  formed by the black dots. The blue region is the feasible space for
  the SDP relaxation~\eqref{sdpprimal}. The orange star is the
  solution to the SDP relaxation, which does not match the global
  solution at the leftmost red square.}
\label{f:threebus_space}
\end{figure*}

Since bus~3 in the three-bus system has zero power injections, it can
be eliminated by adding $R_{13} + \mathbf{j} X_{13}$ and $R_{23} +
\mathbf{j}X_{23}$ to yield an equivalent two-bus system with two
parallel lines.\footnote{Elimination of bus~3 requires that the zero
  power injection at this bus is achieved using an ``open circuit to
  ground''. A ``short circuit to ground'' could also yield zero power
  injections. However, a short circuit at bus~3 results in
  infeasibility of the power flow equations for the loading specified
  in Fig.~\ref{f:threebussystem}. Thus, the feasible space of the
  two-bus system in Fig.~\ref{f:twobussystem} can be directly mapped
  to the feasible space of the three-bus system in
  Fig.~\ref{f:threebussystem}.} The parallel combination of these
lines gives the line impedance $R_{12}^\prime +
\mathbf{j}X_{12}^\prime$ shown in the two-bus system of
Fig.~\ref{f:twobussystem}. Thus, the OPF problems for the two- and
three-bus systems are \emph{equivalent}. The voltage at bus~3 in the
three-bus system can be directly computed from the solution to the
two-bus system. The global solutions are given in
Table~\ref{t:solutions}.

\begin{table}[b]
\vspace{-10pt}
\centering
\caption{Solutions to Two- and Three-Bus Systems (per unit)}
\label{t:solutions}
\begin{tabular}{|c|c|c|}
\hline 
 & Two-Bus System & Three-Bus System\\ \hline\hline
$V_{d1} + \mathbf{j}V_{q1}$ & $1.000 + \mathbf{j}0.000$ & $1.000 + \mathbf{j}0.000$ \\ \hline
$V_{d2} + \mathbf{j}V_{q2}$ & $1.049 - \mathbf{j}0.767$ & $1.049 - \mathbf{j}0.767$\\\hline
$V_{d3} + \mathbf{j}V_{q3}$ & N/A & $0.849 - \mathbf{j}0.586$ \\ \hline
$P_1 + \mathbf{j}Q_1$ & $5.68 - \mathbf{j} 7.77$ & $5.68 - \mathbf{j} 7.77$ \\\hline
$P_2 + \mathbf{j}Q_2$ & $0.0 + \mathbf{j} 12.52$ & $0.0 + \mathbf{j} 12.52$ \\\hline
$P_3 + \mathbf{j}Q_3$ & N/A & $0.0 + \mathbf{j}.0$ \\ \hline
\end{tabular}
\end{table}

The SDP relaxation globally solves the two-bus system. However, for the three-bus system, the relaxation only provides a lower bound that is 22\% less than the true global optimum (i.e., there exists a large relaxation gap). We note that M{\sc atpower}'s interior point solver~\cite{matpower} fails to converge for the three-bus formulation of this problem but successfully solves the two-bus formulation. 

\subsection{Feasible Space Exploration}

Although the OPF problems~\eqref{opf} for the two- and three-bus systems share the same feasible spaces, this is not the case for their SDP relaxations~\eqref{sdpprimal}. This section explores the feasible spaces of these relaxations to illustrate why the SDP relaxation globally solves the two-bus system but fails for the equivalent three-bus system.

Figs.~\ref{f:twobus_space} and~\ref{f:threebus_space} show projections
of the feasible spaces of the two- and three-bus systems,
respectively, in terms of the active power injections. The boundary of
the oval, shown by the black line in Fig.~\ref{f:twobus_space}, is the
feasible space of the OPF problem~\eqref{opf} for varying values of
$P_2$. The region in Fig.~\ref{f:twobus_space} consisting of the oval
and its interior is the feasible space of the SDP relaxation. For the
specified value of $P_2 = 0$, shown by the red dashed line, the OPF
problem has a feasible space consisting of the two red squares at the
intersection of the red dashed line and the black oval. The SDP
relaxation finds the global optimum of~\eqref{opf} (i.e., the leftmost
red square) at the orange star.

In Fig.~\ref{f:threebus_fullspace}, the black dots outline the feasible space of the OPF problem~\eqref{opf} for varying values of $P_2$ and $P_3$, as determined by repeated homotopy calculations~\cite{salam1989}. This feasible space has an ellipsoidal shape with a hole in the interior. The red dashed line corresponds to zero active power injections at buses~2 and~3. The OPF solutions, which are shown by the red squares at the intersection of the exterior of the ellipsoidal shape with the red dashed line, are near the hole in the feasible space. The feasible space of the SDP relaxation, shown by the shaded region, ``stretches over'' this hole in the OPF's feasible space. As seen in Fig.~\ref{f:threebus_zoomspace}, which shows a zoomed view of a cut through $P_3 = 0$, the exterior of the relaxation's feasible space does not match the feasible space of the OPF problem near this hole. Thus, the solution to the SDP relaxation~\eqref{sdpprimal} at the orange star does not match the global solution to the OPF problem~\eqref{opf} at the leftmost red square, and the SDP relaxation is not exact for this formulation. Similar phenomena occur for a range of non-zero active power injections at bus~2.

The hole in the OPF's feasible space is a non-convexity introduced by ``nearby'' problems (i.e., different values of $P_3$) in the three-bus system. Without the additional degrees of freedom associated with bus~3, there is no ``nearby'' non-convexity for the two-bus system. Thus, despite the fact that the OPF problems share the same feasible space (i.e., the red squares in Figs.~\ref{f:twobus_space} and~\ref{f:threebus_space}), the SDP relaxation is exact for the two-bus system but not for the three-bus system.

\section{Moment Relaxations}
\label{l:moment}

By demonstrating that factors other than just physical characteristics determine success or failure of the SDP relaxation, the example in Section~\ref{l:example} motivates the development of tighter convex relaxations that globally solve a broader class of OPF problems. Recognizing that the objective function and all constraints in the OPF problem are polynomial functions of the voltage phasor components enables the application of a hierarchy of convex ``moment'' relaxations from the Lasserre hierarchy for polynomial optimization problems. The moment relaxations, which converge to the global optimum of~\eqref{opf} with increasing relaxation order~\cite{lasserre_book}, generalize the SDP relaxation presented in Section~\ref{l:sdp_relaxation}. This section introduces and illustrates the moment relaxations using the example from Section~\ref{l:example}.

The moment relaxations require definitions beyond those in Section~\ref{l:sdp_relaxation}. Define a vector $x_\gamma$ consisting of all monomials of the voltage components $V_d$ and $V_q$ up to order $\gamma$:
\begin{align}
\nonumber
x_\gamma := & \left[ \begin{array}{ccccccc} 1 & V_{d1} & \ldots & V_{qn} & V_{d1}^2 & V_{d1}V_{d2} & \ldots \end{array} \right. \\ \label{x_d}
& \qquad \left.\begin{array}{cccccc} \ldots & V_{qn}^2 & V_{d1}^3 & V_{d1}^2 V_{d2} & \ldots & V_{qn}^\gamma \end{array}\right]^\intercal.
\end{align}

The moment relaxations are composed of positive semidefinite constraints on \emph{moment} and \emph{localizing} matrices. The symmetric moment matrix $\mathbf{M}_{\gamma}$ is composed of entries $y_\alpha$ corresponding to all monomials $\hat{x}^{\alpha}$ up to order $2\gamma$:
\begin{equation}
\label{eq:real_moment}
\mathbf{M}_\gamma \left\lbrace y \right\rbrace := L_y\left\lbrace x_\gamma^{\vphantom{\intercal}} x_\gamma^\intercal\right\rbrace.
\end{equation}

Symmetric localizing matrices are defined for each constraint of~\eqref{opf}. For a polynomial constraint $g\left(\hat{x}\right) \geqslant 0$ of degree $2\eta$, the localizing matrix is:
\begin{equation}
\label{eq:real_local}
\mathbf{M}_{\gamma - \eta} \left\lbrace g y \right\rbrace := L_y \left\lbrace g x_{\gamma-\eta}^{\vphantom{\intercal}} x_{\gamma-\eta}^{\intercal} \right\rbrace.
\end{equation}

See~\eqref{x2}, \eqref{Moment2}, and \eqref{Local2} for the vector $x_2$, moment matrix $\mathbf{M}_2 \left\lbrace y \right\rbrace$, and the localizing matrix associated with upper voltage magnitude limit $\left(V_2^{\max}\right)^2 - V_{d2}^2 - V_{q2}^2 \geqslant 0$, respectively, for a three-bus OPF problem. Note that the angle reference $V_{q1} = 0$ is used to eliminate $V_{q1}$ in~\eqref{equations2}. These equations use the notation $L_y\left\lbrace V_{d1}^{\alpha_1}V_{d2}^{\alpha_2}V_{d3}^{\alpha_3}V_{q2}^{\alpha_4}V_{q3}^{\alpha_5} \right\rbrace = y^{\alpha_1\alpha_2\alpha_3}_{\hphantom{\alpha_1}\alpha_4\alpha_5}$.

\begin{figure*}[b]
\small
\setcounter{equation}{13}
\hrule
\vspace{10pt}
\begin{subequations}
\label{equations2}
\begin{align} \nonumber \small
x_2 = & \left[\begin{array}{cccccccccccccccc} 1 & V_{d1} & V_{d2} & V_{d3} & V_{q2} & V_{q3} & V_{d1}^2 & V_{d1}V_{d2} & V_{d1}V_{d3} & V_{d1}V_{q2} & V_{d1}V_{q3} & V_{d2}^2 & V_{d2}V_{d3} & V_{d2}V_{q2} & V_{d2}V_{q3} & \ldots \end{array}\right. \\  \label{x2} & \qquad \left.\begin{array}{cccccccccccccccc}\ldots & V_{d3}^2 & V_{d3}V_{q2} & V_{d3}V_{q3} & V_{q2}^2 & V_{q2}V_{q3} & V_{q3}^2\end{array}\right]^\intercal \qquad\; \left[\mathrm{Note\!:}\; \eqref{opf_Vref}\; \mathrm{is\; used\; to\; remove}\; V_{q1} \right]
\end{align}
\vspace{-10pt}
{\arraycolsep=4pt\def\arraystretch{1.25}
\begin{align} \small \label{Moment2} 
& \mathbf{M}_2 \{y \} = L_y\{x_2^{\vphantom{\intercal}} x_2^\intercal\} = \\ 
\nonumber & \left[\begin{array}{ccccccccccccccccccccc}
\tikzmark{left2_0}
y^{000}_{\hphantom{0}00} \tikzmark{right2_0}& \tikzmark{left1top}y^{100}_{\hphantom{0}00} & y^{010}_{\hphantom{0}00} & y^{001}_{\hphantom{0}00} & y^{000}_{\hphantom{0}10} & y^{000}_{\hphantom{0}01}\tikzmark{right1top} & \tikzmark{left2_2top} y^{200}_{\hphantom{0}00} & y^{110}_{\hphantom{0}00} & y^{101}_{\hphantom{0}00} & y^{100}_{\hphantom{0}10} & y^{100}_{\hphantom{0}01} & y^{020}_{\hphantom{0}00} & y^{011}_{\hphantom{0}00} & y^{010}_{\hphantom{0}10} & y^{010}_{\hphantom{0}01} & y^{002}_{\hphantom{0}00} & y^{001}_{\hphantom{0}10} & y^{001}_{\hphantom{0}01} &  y^{000}_{\hphantom{0}20} & y^{000}_{\hphantom{0}11} & y^{000}_{\hphantom{0}02} \tikzmark{right2_2top}\\
\tikzmark{left1left} y^{100}_{\hphantom{0}00} & \tikzmark{left_sdp} y^{200}_{\hphantom{0}00} & y^{110}_{\hphantom{0}00} & y^{101}_{\hphantom{0}00} & y^{100}_{\hphantom{0}10} & y^{100}_{\hphantom{0}01} & y^{300}_{\hphantom{0}00} & y^{110}_{\hphantom{0}00} & y^{201}_{\hphantom{0}00} & y^{200}_{\hphantom{0}10} & y^{200}_{\hphantom{0}01} & y^{120}_{\hphantom{0}00} & y^{111}_{\hphantom{0}00} & y^{110}_{\hphantom{0}10} & y^{110}_{\hphantom{0}01} & y^{102}_{\hphantom{0}00} & y^{101}_{\hphantom{0}10} & y^{101}_{\hphantom{0}01} &  y^{100}_{\hphantom{0}20} & y^{100}_{\hphantom{0}11} & y^{100}_{\hphantom{0}02}\\
y^{010}_{\hphantom{0}00} & y^{110}_{\hphantom{0}00} & y^{020}_{\hphantom{0}00} & y^{011}_{\hphantom{0}00} & y^{010}_{\hphantom{0}10} & y^{010}_{\hphantom{0}01} & y^{210}_{\hphantom{0}00} & y^{120}_{\hphantom{0}00} & y^{111}_{\hphantom{0}00} & y^{110}_{\hphantom{0}10} & y^{110}_{\hphantom{0}01} & y^{030}_{\hphantom{0}00} & y^{021}_{\hphantom{0}00} & y^{020}_{\hphantom{0}10} & y^{020}_{\hphantom{0}01} & y^{012}_{\hphantom{0}00} & y^{011}_{\hphantom{0}10} & y^{011}_{\hphantom{0}01} &  y^{010}_{\hphantom{0}20} & y^{010}_{\hphantom{0}11} & y^{010}_{\hphantom{0}02}\\
y^{001}_{\hphantom{0}00} & y^{101}_{\hphantom{0}00} & y^{011}_{\hphantom{0}00} & y^{002}_{\hphantom{0}00} & y^{001}_{\hphantom{0}10} & y^{001}_{\hphantom{0}01} & y^{201}_{\hphantom{0}00} & y^{111}_{\hphantom{0}00} & y^{102}_{\hphantom{0}00} & y^{101}_{\hphantom{0}10} & y^{101}_{\hphantom{0}01} & y^{021}_{\hphantom{0}00} & y^{012}_{\hphantom{0}00} & y^{011}_{\hphantom{0}10} & y^{011}_{\hphantom{0}01} & y^{003}_{\hphantom{0}00} & y^{002}_{\hphantom{0}10} & y^{002}_{\hphantom{0}01} &  y^{001}_{\hphantom{0}20} & y^{001}_{\hphantom{0}11} & y^{001}_{\hphantom{0}02}\\
y^{000}_{\hphantom{0}10} & y^{100}_{\hphantom{0}10} & y^{010}_{\hphantom{0}10} & y^{001}_{\hphantom{0}10} & y^{000}_{\hphantom{0}20} & y^{000}_{\hphantom{0}11} & y^{200}_{\hphantom{0}10} & y^{110}_{\hphantom{0}10} & y^{101}_{\hphantom{0}10} & y^{100}_{\hphantom{0}20} & y^{100}_{\hphantom{0}11} & y^{020}_{\hphantom{0}10} & y^{011}_{\hphantom{0}10} & y^{010}_{\hphantom{0}20} & y^{010}_{\hphantom{0}11} & y^{002}_{\hphantom{0}10} & y^{001}_{\hphantom{0}20} & y^{001}_{\hphantom{0}11} &  y^{000}_{\hphantom{0}30} & y^{000}_{\hphantom{0}21} & y^{000}_{\hphantom{0}12}\\
y^{000}_{\hphantom{0}01}\tikzmark{right1left} & y^{100}_{\hphantom{0}01} & y^{010}_{\hphantom{0}01} & y^{001}_{\hphantom{0}01} & y^{000}_{\hphantom{0}11} & y^{000}_{\hphantom{0}02} \tikzmark{right_sdp} & y^{200}_{\hphantom{0}01} & y^{110}_{\hphantom{0}01} & y^{101}_{\hphantom{0}01} & y^{100}_{\hphantom{0}11} & y^{100}_{\hphantom{0}02} & y^{020}_{\hphantom{0}01} & y^{011}_{\hphantom{0}01} & y^{010}_{\hphantom{0}11} & y^{010}_{\hphantom{0}02} & y^{002}_{\hphantom{0}01} & y^{001}_{\hphantom{0}11} & y^{001}_{\hphantom{0}02} &  y^{000}_{\hphantom{0}21} & y^{000}_{\hphantom{0}12} & y^{000}_{\hphantom{0}03}\\
\tikzmark{left2_2left} y^{200}_{\hphantom{0}00} & y^{300}_{\hphantom{0}00} & y^{210}_{\hphantom{0}00} & y^{201}_{\hphantom{0}00} & y^{200}_{\hphantom{0}10} & y^{200}_{\hphantom{0}01} & \tikzmark{left2_4} y^{400}_{\hphantom{0}00} & y^{310}_{\hphantom{0}00} & y^{301}_{\hphantom{0}00} & y^{300}_{\hphantom{0}10} & y^{300}_{\hphantom{0}01} & y^{220}_{\hphantom{0}00} & y^{211}_{\hphantom{0}00} & y^{210}_{\hphantom{0}10} & y^{210}_{\hphantom{0}01} & y^{202}_{\hphantom{0}00} & y^{201}_{\hphantom{0}10} & y^{201}_{\hphantom{0}01} &  y^{200}_{\hphantom{0}20} & y^{200}_{\hphantom{0}11} & y^{200}_{\hphantom{0}02}\\
y^{110}_{\hphantom{0}00} & y^{210}_{\hphantom{0}00} & y^{120}_{\hphantom{0}00} & y^{111}_{\hphantom{0}00} & y^{110}_{\hphantom{0}10} & y^{110}_{\hphantom{0}01} & y^{310}_{\hphantom{0}00} & y^{220}_{\hphantom{0}00} & y^{211}_{\hphantom{0}00} & y^{210}_{\hphantom{0}10} & y^{210}_{\hphantom{0}01} & y^{130}_{\hphantom{0}00} & y^{121}_{\hphantom{0}00} & y^{120}_{\hphantom{0}10} & y^{120}_{\hphantom{0}01} & y^{112}_{\hphantom{0}00} & y^{111}_{\hphantom{0}10} & y^{111}_{\hphantom{0}01} &  y^{110}_{\hphantom{0}20} & y^{110}_{\hphantom{0}11} & y^{110}_{\hphantom{0}02}\\
y^{101}_{\hphantom{0}00} & y^{201}_{\hphantom{0}00} & y^{111}_{\hphantom{0}00} & y^{102}_{\hphantom{0}00} & y^{101}_{\hphantom{0}10} & y^{101}_{\hphantom{0}01} & y^{301}_{\hphantom{0}00} & y^{211}_{\hphantom{0}00} & y^{202}_{\hphantom{0}00} & y^{201}_{\hphantom{0}10} & y^{201}_{\hphantom{0}01} & y^{121}_{\hphantom{0}00} & y^{112}_{\hphantom{0}00} & y^{111}_{\hphantom{0}10} & y^{111}_{\hphantom{0}01} & y^{103}_{\hphantom{0}00} & y^{102}_{\hphantom{0}10} & y^{102}_{\hphantom{0}01} &  y^{101}_{\hphantom{0}20} & y^{101}_{\hphantom{0}11} & y^{101}_{\hphantom{0}02}\\
y^{100}_{\hphantom{0}10} & y^{200}_{\hphantom{0}10} & y^{110}_{\hphantom{0}10} & y^{101}_{\hphantom{0}10} & y^{100}_{\hphantom{0}20} & y^{100}_{\hphantom{0}11} & y^{300}_{\hphantom{0}10} & y^{210}_{\hphantom{0}10} & y^{201}_{\hphantom{0}10} & y^{200}_{\hphantom{0}20} & y^{200}_{\hphantom{0}11} & y^{120}_{\hphantom{0}10} & y^{111}_{\hphantom{0}10} & y^{110}_{\hphantom{0}20} & y^{110}_{\hphantom{0}11} & y^{102}_{\hphantom{0}10} & y^{101}_{\hphantom{0}20} & y^{101}_{\hphantom{0}11} &  y^{100}_{\hphantom{0}30} & y^{100}_{\hphantom{0}21} & y^{100}_{\hphantom{0}12}\\
y^{100}_{\hphantom{0}01} & y^{200}_{\hphantom{0}01} & y^{110}_{\hphantom{0}01} & y^{101}_{\hphantom{0}01} & y^{100}_{\hphantom{0}11} & y^{100}_{\hphantom{0}02} & y^{300}_{\hphantom{0}01} & y^{210}_{\hphantom{0}01} & y^{201}_{\hphantom{0}01} & y^{200}_{\hphantom{0}11} & y^{200}_{\hphantom{0}02} & y^{120}_{\hphantom{0}01} & y^{111}_{\hphantom{0}01} & y^{110}_{\hphantom{0}11} & y^{110}_{\hphantom{0}02} & y^{102}_{\hphantom{0}01} & y^{101}_{\hphantom{0}11} & y^{101}_{\hphantom{0}02} &  y^{100}_{\hphantom{0}21} & y^{100}_{\hphantom{0}12} & y^{100}_{\hphantom{0}03}\\
y^{020}_{\hphantom{0}00} & y^{120}_{\hphantom{0}00} & y^{030}_{\hphantom{0}00} & y^{021}_{\hphantom{0}00} & y^{020}_{\hphantom{0}10} & y^{020}_{\hphantom{0}01} & y^{220}_{\hphantom{0}00} & y^{130}_{\hphantom{0}00} & y^{121}_{\hphantom{0}00} & y^{120}_{\hphantom{0}10} & y^{120}_{\hphantom{0}01} & y^{040}_{\hphantom{0}00} & y^{031}_{\hphantom{0}00} & y^{030}_{\hphantom{0}10} & y^{030}_{\hphantom{0}01} & y^{022}_{\hphantom{0}00} & y^{021}_{\hphantom{0}10} & y^{021}_{\hphantom{0}01} &  y^{020}_{\hphantom{0}20} & y^{020}_{\hphantom{0}11} & y^{020}_{\hphantom{0}02}\\
y^{011}_{\hphantom{0}00} & y^{111}_{\hphantom{0}00} & y^{021}_{\hphantom{0}00} & y^{012}_{\hphantom{0}00} & y^{011}_{\hphantom{0}10} & y^{011}_{\hphantom{0}01} & y^{211}_{\hphantom{0}00} & y^{121}_{\hphantom{0}00} & y^{112}_{\hphantom{0}00} & y^{111}_{\hphantom{0}10} & y^{111}_{\hphantom{0}01} & y^{031}_{\hphantom{0}00} & y^{022}_{\hphantom{0}00} & y^{021}_{\hphantom{0}10} & y^{021}_{\hphantom{0}01} & y^{013}_{\hphantom{0}00} & y^{012}_{\hphantom{0}10} & y^{012}_{\hphantom{0}01} &  y^{011}_{\hphantom{0}20} & y^{011}_{\hphantom{0}11} & y^{011}_{\hphantom{0}02}\\
y^{010}_{\hphantom{0}10} & y^{110}_{\hphantom{0}10} & y^{020}_{\hphantom{0}10} & y^{011}_{\hphantom{0}10} & y^{010}_{\hphantom{0}20} & y^{010}_{\hphantom{0}11} & y^{210}_{\hphantom{0}10} & y^{120}_{\hphantom{0}10} & y^{111}_{\hphantom{0}10} & y^{110}_{\hphantom{0}20} & y^{110}_{\hphantom{0}11} & y^{030}_{\hphantom{0}10} & y^{021}_{\hphantom{0}10} & y^{020}_{\hphantom{0}20} & y^{020}_{\hphantom{0}11} & y^{012}_{\hphantom{0}10} & y^{011}_{\hphantom{0}20} & y^{011}_{\hphantom{0}11} &  y^{010}_{\hphantom{0}30} & y^{010}_{\hphantom{0}21} & y^{010}_{\hphantom{0}12}\\
y^{010}_{\hphantom{0}01} & y^{110}_{\hphantom{0}01} & y^{020}_{\hphantom{0}01} & y^{011}_{\hphantom{0}01} & y^{010}_{\hphantom{0}11} & y^{010}_{\hphantom{0}02} & y^{210}_{\hphantom{0}01} & y^{120}_{\hphantom{0}01} & y^{111}_{\hphantom{0}01} & y^{110}_{\hphantom{0}11} & y^{110}_{\hphantom{0}02} & y^{030}_{\hphantom{0}01} & y^{021}_{\hphantom{0}01} & y^{020}_{\hphantom{0}11} & y^{020}_{\hphantom{0}02} & y^{012}_{\hphantom{0}01} & y^{011}_{\hphantom{0}11} & y^{011}_{\hphantom{0}02} &  y^{010}_{\hphantom{0}21} & y^{010}_{\hphantom{0}12} & y^{010}_{\hphantom{0}03}\\
y^{002}_{\hphantom{0}00} & y^{102}_{\hphantom{0}00} & y^{012}_{\hphantom{0}00} & y^{003}_{\hphantom{0}00} & y^{002}_{\hphantom{0}10} & y^{002}_{\hphantom{0}01} & y^{202}_{\hphantom{0}00} & y^{112}_{\hphantom{0}00} & y^{103}_{\hphantom{0}00} & y^{102}_{\hphantom{0}10} & y^{102}_{\hphantom{0}01} & y^{022}_{\hphantom{0}00} & y^{013}_{\hphantom{0}00} & y^{012}_{\hphantom{0}10} & y^{012}_{\hphantom{0}01} & y^{004}_{\hphantom{0}00} & y^{003}_{\hphantom{0}10} & y^{003}_{\hphantom{0}01} &  y^{002}_{\hphantom{0}20} & y^{002}_{\hphantom{0}11} & y^{002}_{\hphantom{0}02}\\
y^{001}_{\hphantom{0}10} & y^{101}_{\hphantom{0}10} & y^{011}_{\hphantom{0}10} & y^{002}_{\hphantom{0}10} & y^{001}_{\hphantom{0}20} & y^{001}_{\hphantom{0}11} & y^{201}_{\hphantom{0}10} & y^{111}_{\hphantom{0}10} & y^{102}_{\hphantom{0}10} & y^{101}_{\hphantom{0}20} & y^{101}_{\hphantom{0}11} & y^{021}_{\hphantom{0}10} & y^{012}_{\hphantom{0}10} & y^{011}_{\hphantom{0}20} & y^{011}_{\hphantom{0}11} & y^{003}_{\hphantom{0}10} & y^{002}_{\hphantom{0}20} & y^{002}_{\hphantom{0}11} &  y^{001}_{\hphantom{0}30} & y^{001}_{\hphantom{0}21} & y^{001}_{\hphantom{0}12}\\
y^{001}_{\hphantom{0}01} & y^{101}_{\hphantom{0}01} & y^{011}_{\hphantom{0}01} & y^{002}_{\hphantom{0}01} & y^{001}_{\hphantom{0}11} & y^{001}_{\hphantom{0}02} & y^{201}_{\hphantom{0}01} & y^{111}_{\hphantom{0}01} & y^{102}_{\hphantom{0}01} & y^{101}_{\hphantom{0}11} & y^{101}_{\hphantom{0}02} & y^{021}_{\hphantom{0}01} & y^{012}_{\hphantom{0}01} & y^{011}_{\hphantom{0}11} & y^{011}_{\hphantom{0}02} & y^{003}_{\hphantom{0}01} & y^{002}_{\hphantom{0}11} & y^{002}_{\hphantom{0}02} &  y^{001}_{\hphantom{0}21} & y^{001}_{\hphantom{0}12} & y^{001}_{\hphantom{0}03}\\
y^{000}_{\hphantom{0}20} & y^{100}_{\hphantom{0}20} & y^{010}_{\hphantom{0}20} & y^{001}_{\hphantom{0}20} & y^{000}_{\hphantom{0}30} & y^{000}_{\hphantom{0}21} & y^{200}_{\hphantom{0}20} & y^{110}_{\hphantom{0}20} & y^{101}_{\hphantom{0}20} & y^{100}_{\hphantom{0}30} & y^{100}_{\hphantom{0}21} & y^{020}_{\hphantom{0}20} & y^{011}_{\hphantom{0}20} & y^{010}_{\hphantom{0}30} & y^{010}_{\hphantom{0}21} & y^{002}_{\hphantom{0}20} & y^{001}_{\hphantom{0}30} & y^{001}_{\hphantom{0}21} &  y^{000}_{\hphantom{0}40} & y^{000}_{\hphantom{0}31} & y^{000}_{\hphantom{0}22}\\
y^{000}_{\hphantom{0}11} & y^{100}_{\hphantom{0}11} & y^{010}_{\hphantom{0}11} & y^{001}_{\hphantom{0}11} & y^{000}_{\hphantom{0}21} & y^{000}_{\hphantom{0}12} & y^{200}_{\hphantom{0}11} & y^{110}_{\hphantom{0}11} & y^{101}_{\hphantom{0}11} & y^{100}_{\hphantom{0}21} & y^{100}_{\hphantom{0}12} & y^{020}_{\hphantom{0}11} & y^{011}_{\hphantom{0}11} & y^{010}_{\hphantom{0}21} & y^{010}_{\hphantom{0}12} & y^{002}_{\hphantom{0}11} & y^{001}_{\hphantom{0}21} & y^{001}_{\hphantom{0}12} &  y^{000}_{\hphantom{0}31} & y^{000}_{\hphantom{0}22} & y^{000}_{\hphantom{0}13}\\
y^{000}_{\hphantom{0}02}\tikzmark{right2_2left} & y^{100}_{\hphantom{0}02} & y^{010}_{\hphantom{0}02} & y^{001}_{\hphantom{0}02} & y^{000}_{\hphantom{0}12} & y^{000}_{\hphantom{0}03} & y^{200}_{\hphantom{0}02} & y^{110}_{\hphantom{0}02} & y^{101}_{\hphantom{0}02} & y^{100}_{\hphantom{0}12} & y^{100}_{\hphantom{0}03} & y^{020}_{\hphantom{0}02} & y^{011}_{\hphantom{0}02} & y^{010}_{\hphantom{0}12} & y^{010}_{\hphantom{0}03} & y^{002}_{\hphantom{0}02} & y^{001}_{\hphantom{0}12} & y^{001}_{\hphantom{0}03} &  y^{000}_{\hphantom{0}22} & y^{000}_{\hphantom{0}13} & y^{000}_{\hphantom{0}04} \tikzmark{right2_4}
\end{array}\right]
\end{align}
\DrawBox[thick, red ]{left1top}{right1top}{}{orange}{dotted}{0.2}{thick}
\DrawBox[thick, red ]{left1left}{right1left}{}{orange}{dotted}{0.2}{thick}
\DrawBox[thick, red ]{left_sdp}{right_sdp}{}{green}{dashed}{0.2}{thick}
\DrawBox[thick, red ]{left2_0}{right2_0}{}{blue}{densely dotted}{0.2}{thick}
\DrawBox[thick, red ]{left2_2top}{right2_2top}{}{blue}{densely dotted}{0.2}{thick}
\DrawBox[thick, red ]{left2_2left}{right2_2left}{}{blue}{densely dotted}{0.2}{thick}
\DrawBox[thick, red ]{left2_4}{right2_4}{}{blue}{densely dotted}{0.2}{thick}
\DrawBox[thick, red ]{left2_0}{right_sdp}{}{}{solid}{0}{ultra thick}
}
{\arraycolsep=4pt\def\arraystretch{1.25}
\begin{align}\nonumber \small & \mathbf{M}_{1}\left\lbrace\left(\left(V_2^{\max}\right)^2 - f_{V2}\right) y \right\rbrace = \\ \label{Local2}  & 
\left(V_2^{\max}\right)^2\!\left[\begin{array}{cccccc}
\tikzmark{left_sdp1}y^{000}_{\hphantom{0}00}\tikzmark{right_sdp1} & y^{100}_{\hphantom{0}00} & y^{010}_{\hphantom{0}00} & y^{001}_{\hphantom{0}00} & y^{000}_{\hphantom{0}10} & y^{000}_{\hphantom{0}01}\\
y^{100}_{\hphantom{0}00} & \tikzmark{left2_1}y^{200}_{\hphantom{0}00} & y^{110}_{\hphantom{0}00} & y^{101}_{\hphantom{0}00} & y^{100}_{\hphantom{0}10} & y^{100}_{\hphantom{0}01}\\
y^{010}_{\hphantom{0}00} & y^{110}_{\hphantom{0}00} & y^{020}_{\hphantom{0}00} & y^{011}_{\hphantom{0}00} & y^{010}_{\hphantom{0}10} & y^{010}_{\hphantom{0}01}\\
y^{001}_{\hphantom{0}00} & y^{101}_{\hphantom{0}00} & y^{011}_{\hphantom{0}00} & y^{002}_{\hphantom{0}00} & y^{001}_{\hphantom{0}10} & y^{001}_{\hphantom{0}01}\\
y^{000}_{\hphantom{0}10} & y^{100}_{\hphantom{0}10} & y^{010}_{\hphantom{0}10} & y^{001}_{\hphantom{0}10} & y^{000}_{\hphantom{0}20} & y^{000}_{\hphantom{0}11}\\
y^{000}_{\hphantom{0}01} & y^{100}_{\hphantom{0}01} & y^{010}_{\hphantom{0}01} & y^{001}_{\hphantom{0}01} & y^{000}_{\hphantom{0}11} & y^{000}_{\hphantom{0}02}\tikzmark{right2_1} \end{array}\right]
\!-\!\left[\begin{array}{cccccc}
\tikzmark{left_sdp2}y^{020}_{\hphantom{0}00}\tikzmark{right_sdp2} & y^{120}_{\hphantom{0}00} & y^{030}_{\hphantom{0}00} & y^{021}_{\hphantom{0}00} & y^{020}_{\hphantom{0}10} & y^{020}_{\hphantom{0}01}\\
y^{120}_{\hphantom{0}00} & \tikzmark{left2_2} y^{220}_{\hphantom{0}00} & y^{130}_{\hphantom{0}00} & y^{121}_{\hphantom{0}00} & y^{120}_{\hphantom{0}10} & y^{120}_{\hphantom{0}01}\\
y^{030}_{\hphantom{0}00} & y^{130}_{\hphantom{0}00} & y^{040}_{\hphantom{0}00} & y^{031}_{\hphantom{0}00} & y^{030}_{\hphantom{0}10} & y^{030}_{\hphantom{0}01}\\
y^{021}_{\hphantom{0}00} & y^{121}_{\hphantom{0}00} & y^{031}_{\hphantom{0}00} & y^{022}_{\hphantom{0}00} & y^{021}_{\hphantom{0}10} & y^{021}_{\hphantom{0}01}\\
y^{020}_{\hphantom{0}10} & y^{120}_{\hphantom{0}10} & y^{030}_{\hphantom{0}10} & y^{021}_{\hphantom{0}10} & y^{020}_{\hphantom{0}20} & y^{020}_{\hphantom{0}11}\\
y^{020}_{\hphantom{0}01} & y^{120}_{\hphantom{0}01} & y^{030}_{\hphantom{0}01} & y^{021}_{\hphantom{0}01} & y^{020}_{\hphantom{0}11} & y^{020}_{\hphantom{0}02}\tikzmark{right2_2} \end{array}\right]
\!-\!\left[\begin{array}{cccccc}
\tikzmark{left_sdp3}y^{000}_{\hphantom{0}20}\tikzmark{right_sdp3} & y^{100}_{\hphantom{0}20} & y^{010}_{\hphantom{0}20} & y^{001}_{\hphantom{0}20} & y^{000}_{\hphantom{0}30} & y^{000}_{\hphantom{0}21}\\
y^{100}_{\hphantom{0}20} & \tikzmark{left2_3}y^{200}_{\hphantom{0}20} & y^{110}_{\hphantom{0}20} & y^{101}_{\hphantom{0}20} & y^{100}_{\hphantom{0}30} & y^{100}_{\hphantom{0}21}\\
y^{010}_{\hphantom{0}20} & y^{110}_{\hphantom{0}20} & y^{020}_{\hphantom{0}20} & y^{011}_{\hphantom{0}20} & y^{010}_{\hphantom{0}30} & y^{010}_{\hphantom{0}21}\\
y^{001}_{\hphantom{0}20} & y^{101}_{\hphantom{0}20} & y^{011}_{\hphantom{0}20} & y^{002}_{\hphantom{0}20} & y^{001}_{\hphantom{0}30} & y^{001}_{\hphantom{0}21}\\
y^{000}_{\hphantom{0}30} & y^{100}_{\hphantom{0}30} & y^{010}_{\hphantom{0}30} & y^{001}_{\hphantom{0}30} & y^{000}_{\hphantom{0}40} & y^{000}_{\hphantom{0}31}\\
y^{000}_{\hphantom{0}21} & y^{100}_{\hphantom{0}21} & y^{010}_{\hphantom{0}21} & y^{001}_{\hphantom{0}21} & y^{000}_{\hphantom{0}31} & y^{000}_{\hphantom{0}22}\tikzmark{right2_3} \end{array}\right]
\end{align}
\DrawBox[thick, red ]{left_sdp1}{right_sdp1}{}{}{solid}{0}{ultra thick}
\DrawBox[thick, red ]{left_sdp2}{right_sdp2}{}{}{solid}{0}{ultra thick}
\DrawBox[thick, red ]{left_sdp3}{right_sdp3}{}{}{solid}{0}{ultra thick}
\DrawBox[thick, red ]{left_sdp1}{right_sdp1}{}{blue}{dotted}{0.2}{thick}
\DrawBox[thick, red ]{left_sdp2}{right_sdp2}{}{blue}{dotted}{0.2}{thick}
\DrawBox[thick, red ]{left_sdp3}{right_sdp3}{}{blue}{dotted}{0.2}{thick}
\DrawBox[thick, red ]{left2_1}{right2_1}{}{blue}{densely dotted}{0.2}{thick}
\DrawBox[thick, red ]{left2_2}{right2_2}{}{blue}{densely dotted}{0.2}{thick}
\DrawBox[thick, red ]{left2_3}{right2_3}{}{blue}{densely dotted}{0.2}{thick}
}
\end{subequations}
\setcounter{equation}{12}
\end{figure*}

The order-$\gamma$ moment relaxation is:
\begin{subequations}
\label{eq:msosr}
\begin{align}
\label{eq:msosr_obj}& \min_{y,\omega} \sum_{k\in\mathcal{G}} \omega_k \qquad \mathrm{subject\; to} \hspace{-150pt} &  \\
\label{eq:msosr_Pmin} &  \mathbf{M}_{\gamma-1}\left\lbrace \left(f_{Pk} - P_k^{\min}\right) y \right\rbrace \succcurlyeq 0 & \forall k\in\mathcal{N} \\
\label{eq:msosr_Pmax} &  \mathbf{M}_{\gamma-1}\left\lbrace \left(P_k^{\max} - f_{Pk} \vphantom{P_k^{\min}}\right) y \right\rbrace \succcurlyeq 0 & \forall k\in\mathcal{N}\\
\label{eq:msosr_Qmin} &  \mathbf{M}_{\gamma-1}\left\lbrace \left(f_{Qk} - Q_k^{\min}\right) y \right\rbrace \succcurlyeq 0 & \forall k\in\mathcal{N}\\
\label{eq:msosr_Qmax} &  \mathbf{M}_{\gamma-1}\left\lbrace \left(Q_k^{\max} - f_{Qk}  \vphantom{P_k^{\min}}\right) y \right\rbrace \succcurlyeq 0 & \forall k\in\mathcal{N}\\
\label{eq:msosr_Vmin} &  \mathbf{M}_{\gamma-1}\left\lbrace \left(f_{Vk} - \left(V_k^{\min}\right)^2\right) y \right\rbrace \succcurlyeq 0 & \forall k\in\mathcal{N}\\
\label{eq:msosr_Vmax} &  \mathbf{M}_{\gamma-1}\left\lbrace \left(\left(V_k^{\max}\right)^2 - f_{Vk}  \vphantom{P_k^{\min}}\right) y \right\rbrace \succcurlyeq 0 & \forall k\in\mathcal{N} \\
\nonumber & \left(1-c_{k1}L_y\left\lbrace f_{Pk} \right\rbrace-c_{k0} + \omega_k \right) \\ \label{eq:msosr_quadcost} & \quad \geqslant \left|\left| \begin{bmatrix} \left(1+c_{k1}L_y\left\lbrace f_{Pk} \right\rbrace+c_{k0}-\omega_k \right) \\ 2\sqrt{c_{k2}}\, L_y\left\lbrace f_{Pk} \right\rbrace \end{bmatrix} \right|\right|_2 \hspace{-22pt} & \forall k \in \mathcal{G}\\
\label{eq:msosr_quadcost_ho} & L_y\left\lbrace f_{Ck} \right\rbrace = \omega_k & \forall k\in\mathcal{G} \\
\label{eq:msosr_Slm_ho} & \mathbf{M}_{\gamma-2}\left\lbrace \left( \left(S_{lm}^{\max}\right)^2 - f_{Slm} \vphantom{P_k^{\min}}\right) y \right\rbrace \succcurlyeq 0 & \forall \left(l,m\right)\in\mathcal{L} \\
\label{eq:msosr_Sml_ho} & \mathbf{M}_{\gamma-2}\left\lbrace \left( \left(S_{lm}^{\max}\right)^2 - f_{Sml} \vphantom{P_k^{\min}}\right) y \right\rbrace \succcurlyeq 0 & \forall \left(l,m\right)\in\mathcal{L} \\
\label{eq:msosr_Slm} & S_{lm}^\mathrm{max} \geqslant \left|\left| \begin{bmatrix} L_y\left\lbrace f_{Plm} \right\rbrace \\ L_y\left\lbrace f_{Qlm} \right\rbrace  \end{bmatrix} \right|\right|_2 & \forall \left(l,m \right) \in \mathcal{L} \\
\label{eq:msosr_Sml} & S_{lm}^\mathrm{max} \geqslant \left|\left| \begin{bmatrix} L_y\left\lbrace f_{Pml} \right\rbrace \\ L_y\left\lbrace f_{Qml} \right\rbrace \end{bmatrix} \right|\right|_2 & \forall \left(l,m \right) \in \mathcal{L}  \\
\label{eq:msosr_Msdp} & \mathbf{M}_\gamma \{y\} \succcurlyeq 0 & \\
\label{eq:msosr_y0} & y_{00\ldots 0} = 1 & \\
\label{eq:msosr_Vref} & y_{\star\star\ldots\star\rho\star\ldots\star} = 0 & \hspace{-20pt}\rho = 1,\ldots,2\gamma.
\end{align}
\end{subequations}

\noindent where $\rho$ in the angle reference constraint~\eqref{eq:msosr_Vref} is in the index $n+1$, which corresponds to the variable $V_{q1}$. In the same way as~\eqref{sdpprimal}, the angle reference can alternatively be used to eliminate all terms corresponding to $V_{q1}$. 

As for the SDP relaxation, the globally optimal voltage phasors can be extracted using~\eqref{Vstar} from a solution to~\eqref{eq:msosr} that satisfies the condition $\mathrm{rank}\left(L_y\left\lbrace \hat{x}\hat{x}^\intercal \right\rbrace\right) = 1$.

The order $\gamma$ of the moment relaxation~\eqref{eq:msosr} must be greater than or equal to half of the degree of any polynomial in the OPF problem~\eqref{opf}. This suggests that $\gamma \geqslant 2$ due to the fourth-order polynomials resulting from the objective function~\eqref{opf_obj} and the apparent power line flow constraints~\eqref{opf_Slm} and~\eqref{opf_Sml}. However, as in the SDP relaxation~\eqref{sdpprimal}, these can be rewritten using a Schur complement~\cite{andersen2014} to allow $\gamma \geqslant 1$. Experience suggests that implementing~\eqref{opf_obj}, \eqref{opf_Slm}, and~\eqref{opf_Sml} both directly and with a Schur complement formulation, as shown in~\eqref{eq:msosr_quadcost} and~\eqref{eq:msosr_quadcost_ho} for the quadratic objective function and~\eqref{eq:msosr_Slm_ho}--\eqref{eq:msosr_Sml}, gives superior results for $\gamma \geqslant 2$. (Constraints~\eqref{eq:msosr_quadcost_ho}, \eqref{eq:msosr_Slm_ho}, and~\eqref{eq:msosr_Sml_ho} are not enforced for $\gamma = 1$.)

The second-order relaxation's moment matrix $\mathbf{M}_2 \{y\}$ is shown in~\eqref{Moment2}. The upper limit on the voltage magnitude at bus~2 in~\eqref{eq:msosr_Vmax} corresponds to a positive semidefinite constraint on the localizing matrix shown in~\eqref{Local2}.

\begin{figure*}[t]
\centering
\subfloat[\hspace{0pt}Projection of the $2^{\mathrm{nd}}$-Order Moment Relaxation's Feasible Space]{\includegraphics[trim=0.2cm 0cm 0.5cm 5cm, clip=true,totalheight=0.20\textheight]{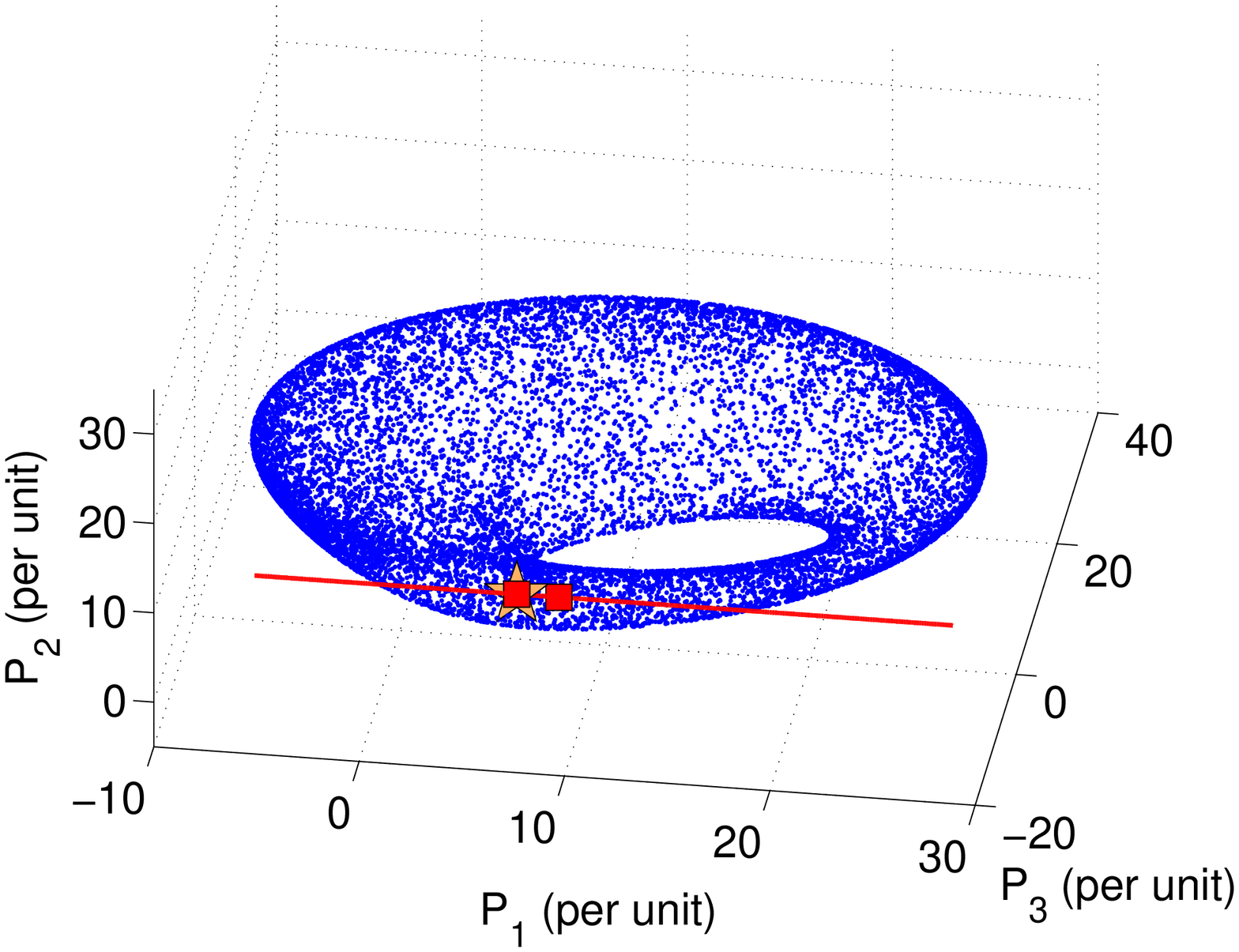}\label{f:threebus_fullspace_msdp2}}\hspace{15pt}
\subfloat[\hspace{0pt}Projection of the Feasible Space shown in Fig.~\ref{f:threebus_fullspace_msdp2} with $P_{3} = 0$]{\includegraphics[totalheight=0.24\textheight]{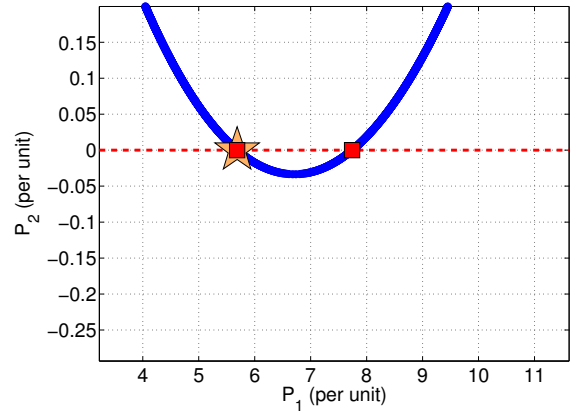}\label{f:threebus_zoomspace_msdp2}}
\caption{Projection of the Second-Order Moment Relaxation's Feasible Space for the Three-Bus System. The feasible space for the OPF problem~\eqref{opf} is denoted by the red squares. The second-order moment relaxation gives the global optimum at the orange star. The second-order moment relaxation was exact for all scenarios tested.}
\label{f:relaxations_space_sdp}
\end{figure*}

Fig.~\ref{f:relaxations_space_sdp} shows a projection of the feasible
space in terms of active power injections for the second-order moment
relaxation of the three-bus system in
Section~\ref{l:example}. The points in this figure were
  obtained by gridding the $P_1$--$P_2$--$P_3$ space, and associating
  with each grid point a quadratic objective function that achieved
  its minimum at that point. The relaxation \eqref{eq:msosr} was
  solved for each of those objective functions while allowing the
  loading conditions to vary (i.e., the constraints on $P_2$ and $P_3$
  were released). The second-order moment relaxation globally solved
all these scenarios, with the resulting feasible space in
Fig.~\ref{f:relaxations_space_sdp} seemingly equivalent to the space
illustrated by the black dots in Fig.~\ref{f:threebus_space}. Since
the power injections result from a non-linear transformation of the
voltage components given by the power flow
equations~\eqref{opf_balance}, the second-order moment relaxation can
represent the non-convex space of power injections while maintaining
convexity in the decision variables $y_{\alpha}$.

All polynomials in the OPF problem have only even-order monomials
(i.e., $\hat{x}^\alpha$ such that $\left|\alpha\right|$ is even, where
$\left|\,\cdot\,\right|$ indicates the one-norm). Odd-order terms in
the moment relaxations are therefore unnecessary: all $y_\alpha$ such
that $\left|\alpha\right|$ is odd can be set to zero without violating
any constraints or changing the objective value. For instance, the
positive semidefinite constraint on the second-order relaxation's
moment matrix, $\mathbf{M}_2 \{y\} \succcurlyeq 0$, is
equivalent to positive semidefinite constraints on two submatrices:
the diagonal block corresponding to the degree-two monomials (i.e.,
$\left|\alpha\right| = 2$), which is identified by the green dashed
highlighting in~\eqref{Moment2}, and the terms corresponding to the
degree-zero, off-diagonal degree-two, and degree-four monomials (i.e., $\left|\alpha\right| = 2k$ for some $k \in \mathbb{N}$), which are identified by the blue
dotted highlighting in~\eqref{Moment2}.

The first-order localizing ``matrices'' corresponding to the
constraints~\eqref{sdp_P}--\eqref{sdp_V} are, in fact,
scalars.\footnote{Observe that $L_y\left\lbrace g\left(\hat{x}\right)
  x_0^{\vphantom{\intercal}} x_0^\intercal \right\rbrace =
  L_y\left\lbrace g\left(\hat{x}\right) \right\rbrace$ since $x_0 =
  1$.} The corresponding scalar constraints in the first-order
relaxation~\eqref{eq:msosr_Pmin}--\eqref{eq:msosr_Vmax} are equivalent
to the linear constraints in the SDP
relaxation~\eqref{sdp_P}--\eqref{sdp_V}. The moment matrix in the
first-order relaxation has all terms $y_{\alpha}$ such that
$\left|\alpha\right| \leqslant 2$ (the diagonal block surrounded by
the black line in~\eqref{Moment2}), whereas the SDP relaxation has all
terms $y_{\alpha}$ such that $\left|\alpha\right| = 2$ (the diagonal
block with green dashed highlighting
in~\eqref{Moment2}). The degree-one terms (the terms with
  orange highlighting in~\eqref{Moment2}) have odd $|\alpha|$ and so
  are unnecessary, as discussed earlier. Therefore, since $y_{00\ldots
    0} \geqslant 0$ by~\eqref{eq:msosr_y0}, the positive semidefinite
constraint on the first-order relaxation's moment
matrix~\eqref{eq:msosr_Msdp} is equivalent to the positive
semidefinite constraint in the SDP relaxation~\eqref{sdp_W}. With an
equivalent feasible space and objective function, the SDP relaxation
in~\eqref{sdpprimal} is the same as the first-order $(\gamma = 1)$
moment relaxation \eqref{eq:msosr}.


\begin{figure*}[t]
\centering
\subfloat[{\hspace{0pt}Projection of the $2^{\mathrm{nd}}$-Order Mixed SDP/SOCP Relaxation's Feasible Space}]{\includegraphics[trim=0.2cm 0cm 0.5cm 5cm, clip=true,totalheight=0.20\textheight]{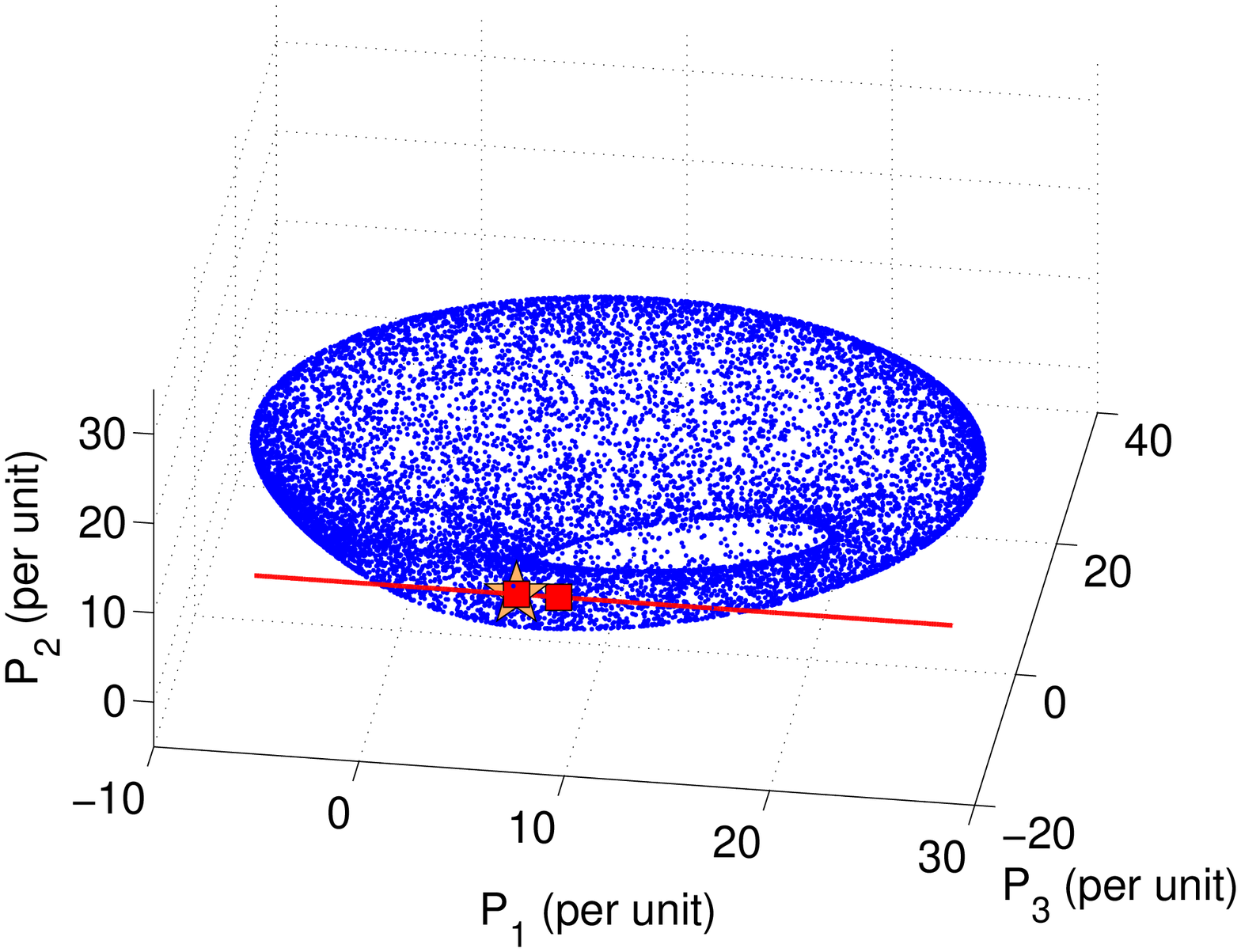}\label{f:threebus_fullspace_sdpsocp2}}
\hspace{15pt}
\subfloat[\hspace{0pt}Projection of the Feasible Space shown in Fig.~\ref{f:threebus_fullspace_sdpsocp2} with $P_{3} = 0$]{\includegraphics[totalheight=0.24\textheight]{threebus_zoomspace_cut_order2.eps}\label{f:threebus_zoomspace_sdpsocp2}}
\caption{Projection of the Second-Order Mixed SDP/SOCP Relaxation's Feasible Spaces for the Three-Bus System. The feasible space for the OPF problem~\eqref{opf} is denoted by the red squares. The second-order mixed SDP relaxation gives the global optimum at the orange star. Fig.~\ref{f:threebus_zoomspace_sdpsocp2} shows that the second-order mixed SDP/SOCP relaxation is exact for the points near the specified scenario. However, this was not the case for all scenarios: Fig.~\ref{f:threebus_fullspace_sdpsocp2} shows that the second-order mixed SDP/SOCP relaxation includes some points in the ``hole'' in the feasible space for which $\mathrm{rank}\left(L_y\left\lbrace \hat{x}\hat{x}^\intercal\right\rbrace\right) > 1$.}
\label{f:relaxations_space_sdpsocp}
\end{figure*}

The moment matrix for the lower-order relaxation $\mathbf{M}_{\gamma-1} \{y\}$ is contained in the upper-left diagonal block of $\mathbf{M}_{\gamma} \{y\}$. Likewise, the upper-left diagonal block of the higher-order localizing matrices contain the lower-order localizing matrices. (The first-order matrices are contained within the solid black outlines in the second-order matrices in \eqref{Moment2} and~\eqref{Local2}.) Since a necessary condition for a matrix to be positive semidefinite is positive semidefiniteness of all principal submatrices, the moment relaxations form a \emph{hierarchy} where higher-order constraints imply the lower-order constraints.

Adding a rank-constraint $\mathrm{rank}\left(L_y\left\lbrace \hat{x}\hat{x}^\intercal\right\rbrace\right) = 1$ to the SDP relaxation~\eqref{sdpprimal} yields a non-convex problem equivalent to the OPF problem~\eqref{opf}. The SDP formulation~\eqref{sdpprimal} can thus be understood in terms of a rank relaxation. The higher-order moment relaxations generalize this approach by introducing constraints that are redundant in the OPF problem~\eqref{opf} but strengthen the moment relaxations. Consider $g\left(\hat{x}\right) x_{\gamma-\eta}^{\vphantom{\intercal}} x_{\gamma-\eta}^\intercal \succcurlyeq 0$, where $g\left(\hat{x}\right) \geqslant 0$ is a generic constraint in the OPF problem~\eqref{opf} with degree $2\eta$. The rank-one matrix $x_{\gamma-\eta}^{\vphantom{\intercal}} x_{\gamma-\eta}^\intercal$ is positive semidefinite by construction, and the scalar constraint $g\left(\hat{x}\right)$ is non-negative. Thus, their product is a rank-one positive semidefinite matrix. Relaxing to $L_y\left\lbrace g\left(\hat{x}\right) x_{\gamma-\eta}^{\vphantom{\intercal}} x_{\gamma-\eta}^\intercal\right\rbrace \succcurlyeq 0$ (i.e., eliminating the rank constraint implied by $x_{\gamma-\eta}^{\vphantom{\intercal}} x_{\gamma-\eta}^\intercal$) results in the localizing matrix constraint.

The computational difficulty of solving the moment relaxations grows quickly with the relaxation order due to the size of the positive semidefinite matrix constraints. After elimination of $V_{q1}$ using the angle reference constraint, the size of the moment matrix~\eqref{eq:msosr_Msdp} for the order-$\gamma$ relaxation of a $n$-bus system is $\left(2n-1+\gamma\right)! / \left( \left(2n-1\right)! \gamma!\right)$. For instance, the third-order relaxation of a 10 bus system has a matrix with size $1,\!540\times 1,\!540$. The ``dense'' formulation of the second-order relaxation is limited to solving problems with less than approximately ten buses~\cite{cedric,pscc2014,ibm_paper}. By exploiting sparisty using techniques analogous to those for the SDP relaxation~\cite{waki2006}, the second-order relaxation is computationally tractable for systems with up to approximately 40 buses~\cite{molzahn_hiskens-sparse_moment_opf}. Extension to larger systems is possible by both exploiting sparsity and only applying the computationally intensive higher-order constraints to specific ``problematic'' buses~\cite{molzahn_hiskens-sparse_moment_opf,cdc2015}.

\section{Mixed SDP/SOCP Relaxation Hierarchy}
\label{l:mixed_sdpsocp}

The moment relaxations globally solve many OPF problems but are computationally challenging. First proposed in~\cite{powertech2015}, a ``mixed SDP/SOCP'' hierarchy is tighter than the first-order moment relaxation but more tractable than the higher-order relaxations. This section describes this mixed SDP/SOCP hierarchy in the context of the example problem from Section~\ref{l:example}.

\setcounter{equation}{14}

The mixed SDP/SOCP hierarchy further relaxes the SDP constraints in
the higher-order moment relaxations~\eqref{eq:msosr} to less stringent
SOCP constraints. To ensure that the mixed SDP/SOCP relaxations are at
least as tight as the first-order moment relaxation, positive
semidefinite constraints are enforced for the diagonal block of the
moment matrix that corresponds to degree-two monomials (i.e.,
$y_\alpha$ such that $\left|\alpha\right| = 2$, which are contained
within the diagonal block highlighted in green dashed lines
in~\eqref{Moment2}.) We relax the higher-order constraints using a
necessary condition for a matrix to be positive
semidefinite. Specifically, a necessary but not sufficient condition
for $\mathbf{W} \succcurlyeq 0$ is given by the constraints:
\begin{subequations}
\label{socp_necessary}
\begin{align}
& \mathbf{W}_{ii} \geqslant 0 & i = 1,\ldots,2n, \\ 
& \mathbf{W}_{ii} \mathbf{W}_{kk} \geqslant \left|\mathbf{W}_{ik}\right|^2 & \forall \left\lbrace \left(i,k\right) \mid k > i\right\rbrace.
\end{align}
\end{subequations}


The mixed SDP/SOCP hierarchy enforces the higher-order constraints in the moment and localizing matrices using~\eqref{socp_necessary} for~\eqref{eq:msosr_Pmin}--\eqref{eq:msosr_Vmax}, \eqref{eq:msosr_Slm_ho}--\eqref{eq:msosr_Sml_ho}, and \eqref{eq:msosr_Msdp}. Since the terms corresponding to odd-order monomials can be set to zero, this reduces to enforcing the SOCP constraints on the submatrices corresponding to those highlighted in blue in~\eqref{Moment2} and~\eqref{Local2} for the three-bus system. 


Since SOCP constraints have significant computational advantages over
SDP constraints, the mixed SDP/SOCP relaxation is more tractable than
the formulation of the moment relaxations given in
Section~\ref{l:moment}. Further, it is only necessary to enforce the
SOCP constraints that correspond to terms in the higher-order matrices
that appear in a localizing matrix constraint. This provides
additional computational advantages when combined with the approach of
selectively applying the higher-order relaxation
constraints~\cite{molzahn_hiskens-sparse_moment_opf}. See~\cite{powertech2015}
for detailed numerical results demonstrating speed increases between a
factor of 1.13 and 18.70 compared to the moment relaxations.

Fig.~\ref{f:relaxations_space_sdpsocp} shows the feasible space of
power injections for the second-order mixed SDP/SOCP
relaxation. This figure was produced using the same
  gridding procedure employed in
  Fig.~\ref{f:relaxations_space_sdp}. The relaxation is exact for the
specific loading condition $P_2 = P_3 = 0$ considered in
Section~\ref{l:example} and for nearby loading conditions (see the
zoomed-in view of the feasible space shown in
Fig.~\ref{f:threebus_zoomspace_sdpsocp2}). However, in contrast to the
moment relaxations implemented with SDP constraints alone, illustrated
in Fig.~\ref{f:relaxations_space_sdp}, the mixed SDP/SOCP relaxation
was not exact for all scenarios. This is evident by the points that
lie in the ``hole'' in the feasible space (i.e., the
points in Fig.~\ref{f:threebus_fullspace_sdpsocp2} that are not in
Fig.~\ref{f:threebus_fullspace_msdp2}).\footnote{The feasible spaces
  for the third- and fourth-order mixed SDP/SOCP relaxations also had
  points in this hole.} As expected, mixed SDP/SOCP
  relaxations are generally not as tight as the moment relaxations
  in~\eqref{eq:msosr} which use only SDP constraints.


\section{Conclusion}
\label{l:conclusion}

An SDP relaxation globally solves many OPF problems which do not satisfy any existing sufficient conditions that assure exactness of the relaxation. This motivates the development of broader sufficient conditions, with a common conjecture being that some physical characteristics of the OPF problem can determine success or failure of the SDP relaxation. This paper has presented a small example OPF problem with two equivalent formulations. The SDP relaxation globally solves only one of the two formulations. This suggests that strictly physically based sufficient conditions for exactness of the
SDP relaxation of the OPF problem cannot predict the relaxation's success or failure for all OPF problems.

The inability to develop universal, physically based sufficient conditions for success of the SDP relaxation motivates researching more sophisticated convex relaxations. We use the small example problem to illustrate two recently developed convex relaxation hierarchies: ``moment'' relaxations from the Lasserre hierarchy for polynomial optimization and a mixed SDP/SOCP hierarchy derived by relaxing the higher-order constraints in the moment relaxations. Both of these hierarchies generalize the SDP relaxation in order to enable global solution of a broader class of OPF problems.

\bibliographystyle{IEEEtran}
\bibliography{illustrative_example}

\newpage
\begin{IEEEbiography}[{\includegraphics[width=1in,height=1.25in,clip,keepaspectratio]{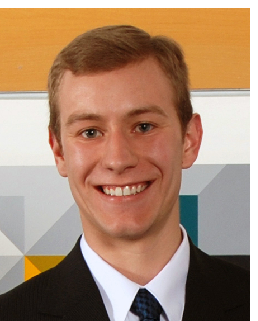}}]{Daniel K. Molzahn}
(S'09-M'13) is a Computational Engineer at Argonne National Laboratory. He recently completed the Dow Postdoctoral Fellow in Sustainability at the University of Michigan, Ann Arbor. He received the B.S., M.S., and Ph.D. degrees in electrical engineering and the Masters of Public Affairs degree from the University of Wisconsin--Madison, where he was a National Science Foundation Graduate Research Fellow. His research focuses on power system optimization and control.
\end{IEEEbiography}
\begin{IEEEbiography}[{\includegraphics[width=1in,height=1.25in,clip,keepaspectratio]{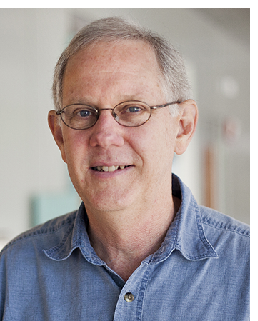}}]{Ian A. Hiskens}
(F'06) is the Vennema Professor of Engineering with the Department of Electrical Engineering and Computer Science, University of Michigan, Ann Arbor. He has held prior appointments in the electricity supply industry (for ten years) and various universities in Australia and the United States. His research focuses on power system analysis, in particular the modelling, dynamics, and control of large-scale, networked, nonlinear systems. His recent activities include integration of renewable generation and new forms of load.

Prof. Hiskens is actively involved in various IEEE societies, and is VP-Finance of the IEEE System Council. He is a Fellow of Engineers Australia and a Chartered Professional Engineer in Australia.
\end{IEEEbiography}

\vfill

\end{document}